\documentclass[11pt]{elsart-maktor}
\usepackage{amsmath}
\renewcommand{\and}{et}
\usepackage{amssymb,
amsfonts,amssymb,euscript,mathrsfs,pifont
}
\usepackage[english,francais]{babel}
\usepackage{amsxtra}
\usepackage{amscd}
\usepackage[T1]{fontenc}                
\usepackage[all]{xy}
\usepackage{stmaryrd}

\usepackage[utf8]{inputenc}
\usepackage[T1]{fontenc}

\usepackage{amssymb,
amsfonts,amssymb,euscript,mathrsfs,pifont
}
\usepackage[english,francais]{babel}
\usepackage{amsxtra}
\usepackage{amscd}
\usepackage[T1]{fontenc}                
\usepackage[all]{xy}

 \usepackage{stmaryrd}

\renewcommand{\and}{et}
\usepackage{amssymb,
amsfonts,amssymb,euscript,mathrsfs,pifont
}
\usepackage[english,francais]{babel}
\usepackage{amsxtra}
\usepackage{amscd}
\usepackage[T1]{fontenc}                
\usepackage[all]{xy}

 \usepackage{stmaryrd}

 \DeclareMathAlphabet{\got}{U}{euf}{m}{n}     
\DeclareMathAlphabet{\mat}{U}{msb}{m}{n}     
\DeclareMathAlphabet{\mathbold}{OML}{cmm}{bx}{it} 

\makeatletter

\renewcommand{\section}{\@startsection{section}{1}%
\z@{.5ex \@plus.7ex}{-.5em}%
{\normalfont\bfseries}}

\renewcommand{\subsection}{\@startsection{subsection}{2}%
\z@{.5ex \@plus.7ex}{-.5em}%
{\normalfont\bfseries}}

\renewcommand{\subsubsection}{\@startsection{subsubsection}{3}%
\z@{.5ex\@plus.7ex}{-.5em}%
{\normalfont\itshape}}
\makeatother
\numberwithin{equation}{section}

\newtheorem{tho}{Theorem}[subsection]
\newtheorem{pro}{Proposition}[subsection]
\newtheorem{lem}{Lemma}[subsection]
\newtheorem{defi}{Definition}[subsection]
\newtheorem{co}{Corollary}[subsection]
\newtheorem{remm}{Remarks}[subsection]

\newcommand{\goth}{\mathfrak}

{\par\rm} 

\newenvironment{theo*}{\vskip 1em\bf Th\'eor\`eme.\it}{\par\rm}
\newenvironment{dem}{\vskip 1em{\it D\'emonstration} :}%
{\unskip\hfill\null\nobreak\hfill\carre\vskip1em\par}

\newenvironment{Proof}{\vskip 1em{\it Proof}.}%
{\unskip\hfill\null\nobreak\hfill\carre\vskip1em\par}

\newcommand{\carre}{\rule{1ex}{1ex}} 

\renewcommand{\thefootnote}{\fnsymbol{footnote}}

\makeatletter
\providecommand*{\diff}%
{\@ifnextchar^{\DIfF}{\DIfF^{}}}
\def\DIfF^#1{%
\mathop{\mathrm{\mathstrut d}}%
\nolimits^{#1}\gobblespace}
\def\gobblespace{%
\futurelet\diffarg\opspace}
\def\opspace{%
\let\DiffSpace\!%
\ifx\diffarg(%
\let\DiffSpace\relax
\else
\ifx\diffarg[%
\let\DiffSpace\relax
\else
\ifx\diffarg\{%
\let\DiffSpace\relax
\fi\fi\fi\DiffSpace}

\providecommand*{\eu}%
{\ensuremath{\mathrm{e}}}

\providecommand*{\iu}%
{\ensuremath{\mathrm{i}}}

\newcommand{\Ind}{\mathop{\mathrm{Ind}}\nolimits}

\newcommand{\vol}{\mathop{\mathrm{vol}}}

\def \cit { \mathbb C}

\def \N {\mathbb N}  
\def \Z {\mathbb Z }

          \def\var{\varphi}
          \def \b{\goth b}
\def\g{ \goth g}                  
                       
           \def \h{\goth h}
    \def \‡{\varsigma}      
            \def \O { \mathcal O}

\def \O{\mathcal O}      
\def\var{\varphi}      
    \def \a{\goth a}

       \def\H{\mathcal H}
           
     \def \e { \varepsilon }

\begin{document}

\begin{frontmatter}
\title{ Harmonic Analysis on  $p$-adic Unipotent Groups}

\author{Khemais Maktouf}
\address{University of Sousse, Laboratory of  Mathematics and Physics, Special Functions and Applications (LR 11 ES 35), ESSTHS 4002 Sousse, 
 Tunisia}
\ead{khemais.maktouf@fsm.rnu.tn}

\selectlanguage{french}

\selectlanguage{english}

\begin{abstract}

Let $G$ be a $p$-adic unipotent group and let $H$ be an algebraic subgroup of $G$. Given a unitary character $\chi$ of $H$, we consider the monomial representation $\operatorname{Ind}_{H}^{G}\chi$.
 To the character $\chi$, we associate the affine subspace
$$
X_\chi=f+\mathfrak h^\perp\subset\mathfrak g^*$$
where $f$ is a linear form  associated with $\chi$. We describe the irreducible representations occurring in the central decomposition
 of $\operatorname{Ind}_{H}^{G}\chi$  
in terms of the coadjoint orbits in $\mathfrak g^*$ that meet $X_\chi$. We also determine the corresponding measure class on the unitary dual $\widehat G$. Furthermore, we study   the multiplicities of the irreducible representations in this decomposition. As an application,  we compute the Plancherel formula for $\operatorname{Ind}_{H}^{G}\chi$ when $H$ is a normal algebraic  subgroup of $G$ and the multiplicities are finite. 

\end{abstract}
\end{frontmatter}
\selectlanguage{french}

\section{Introduction :}

 Let $\bf {G}$ be a unipotent algebraic group defined over    a $p$-adic field  $F$ of characteristic zero. The set of rational points $G$ of $\bf {G}$  is a unipotent $p$-adic group  with Lie algebra  $\g$ and the exponential map $\exp$ realizes   a diffeomorphism from  $\g$ onto  $G$. 
The group $G$ acts on its Lie algebra by the adjoint representation and on  the  dual space $\g^*$ by the coadjoint representation. 
According to the orbit method,  the unitary dual $\widehat{G}$ of $G$  is parametrized by 
 $G \backslash \g^*$,  the space of  coadjoint orbits. Thus for $\rho \in \widehat{G}$, we denote by $\Omega_\rho$ its corresponding coadjoint orbit.

Let $H$  be an algebraic subgroup of $G$,  and let $\chi$  be a unitary character of $H$. A natural class of unitary representations of $G$ is given by the monomial representations
$
\operatorname{Ind}_{H}^{G}\chi$.

Let 
$\h$ be  the Lie algebra of $H$. Then there exists  a linear form
$f\in\mathfrak{g}^{*}$ such that $\chi=\chi_{f}$, where
\begin{equation*}
\chi_{f}(\exp X)=\psi_F(\langle f,X\rangle)\mbox{, }X\in\mathfrak{h},
\end{equation*}
 and $\psi_F$ is a fixed   nontrivial continuous  character of $(F,+)$ into $\cit\setminus\{0\}$.
The Lie algebra $\mathfrak{h}$ must be subordinate to $f$, that is,
$f([\mathfrak{h},\mathfrak{h}])=0$.

When $\h$ is  a polarization of the linear form $f$, the induced representation is irreducible and corresponds, under the orbit method, to the coadjoint orbit $G \cdot f$.  In general,  however, if $\h$ is not a polarization of  $f$, the representation $\operatorname{Ind}_{H}^{G}\chi$ need not be irreducible and its decomposition into irreducible unitary representations becomes a natural problem.

The purpose of this paper is to study the central decomposition of monomial representation  $\operatorname{Ind}_{H}^{G}\chi$ of the group $G$. Since such group is of  type I, every unitary representation admits a central decomposition, which in this setting can be described as a direct integral of irreducible representations.

We investigate,  in term of coadjoint orbits meeting the affine subspace  $X_\chi = f + \h^\perp$, consisting of linear forms whose restriction to the Lie algebra $\h$ corresponds to $f$,   which  irreducible representations occur in the central decomposition of $\operatorname{Ind}_{H}^{G}\chi$. We also determine the measure class of this central desintegration and study   the multiplicities of the irreducible representations in it.  More precisely,  
let $ \mu_{\h} $ be a Haar measure on $  \h^\perp$,  and let $ \mu_{f, \h}$ 
denote the Borel measure on $X_\chi$, obtained by   translating   $\mu_{\h}$ by $f$. 
We regard this measure naturally as a Borel measure on $\g^*$.
Let $ \nu_{f, \h}$  be a  finite Borel measure on $X_\chi$ that is equivalent to 
$ \mu_{f, \h}$. We obtain : 
\begin{equation}\label{princ} 
   \operatorname{Ind}_{H}^{G}\chi = \int^\oplus_{G.X_\chi/G}n_\Omega \pi_\Omega  \diff\dot\nu_{f, \h}^G(\Omega)
	\end{equation}
where $G.X_\chi/G$ is the set of $G$-coadjoint orbits that meet $X_\chi$, $\dot\nu_{f, \h}^G$ is the  push-forward of    $\nu_{f, \h}$ under $X_\chi\longrightarrow G.X_\chi/G$, and $n_\Omega \in \N \cup \{+\infty\} $.

Let  $\Xi$ denote the set of elements $l \in   f + \h^\perp$ such that  $\h + \g(l)$ is a maximal   isotropic subspace (for $\beta_l$).

To compute the multiplicities $n_\Omega$ appearing in formula \ref{princ}, we prove the existence of a nonempty Zariski-open subset $\mathcal U$ of $f + \h^\perp$   that is  a union of intersections with $G$-orbits, and   satisfies the following property : for every   $l \in \mathcal U$,  
 $G.l \cap (f + \h^\perp) $  is a smooth variety of dimension equal to     $\dim(\g.l \cap \h^\perp)$. Moreover 
$\Xi \cap \mathcal U$ is  a Zariski open subset  in $f + \h^\perp$.
We then obtain : 
\begin{itemize}
	\item if $\Xi \cap \mathcal U \neq \varnothing$, then $$ n_\Omega= cad\{\mbox{$H$-orbits in $\Omega \cap X_\chi$}\}, \, \, \forall \Omega \cap \mathcal U \neq \varnothing ; $$ 
	\item  if $\Xi \cap \mathcal U  =  \varnothing$, then $$ n_\Omega= +\infty, \forall\Omega.$$ 
\end{itemize}

As an application,  we derive the Plancherel formula for $\operatorname{Ind}_{H}^{G}\chi$ when $H$ is a normal subgroup of $G$ and the multiplicities are finite. To this end,  we study the notion of distribution vectors associated with smooth representations of $G$.


In the real case,  these results are  established  by Corwin, Greenleaf, and Grelaud \cite{corwin-greenleaf-grelaud-1987} and by Lipsman \cite{lipsman-1989} for connected, simply connected nilpotent
Lie group. They  determine the multiplicities of  irreducible representations
 in the central desintegration of $\Ind_{H}^{G}\chi_f$. 

\section{ Notations and generalities:  }\label{notation}
\subsection{}\label{intr} 
Along this paper, $F$ denotes a nonarchimedean local field of zero
characteristic and residual characteristic $p$. Let $\mathcal{O}$
be the ring of integers of $F$, $\mathcal{P}$ its maximal ideal, and
$\varpi$ a uniformizer  of $F$, that is, 
$\mathcal{P}=\varpi\mathcal{O}$. Then, the residual class field
$\mathcal{O}/\mathcal{P}=\mathbb{F}_{q}$ is a finite field of
characteristic $p$ and cardinal $q$, a power of $p$.   Let $v :
k^{\times}\rightarrow\mathbb{Z}$ be the discrete valuation such that
$v(\varpi)=1$.  We extend $v$ to $0$ by
putting $v(0)=+\infty$.

We fix a nontrivial unitary character $\psi_F$ of the additive group $F$  and denote by  $\lambda_\psi $ its  conductor, that is,   $\lambda_\psi $ is the smallest integer such that 
\begin{equation*}
\psi_{|\goth p^{\lambda_\psi}}\equiv 1 \mbox{ and  } \psi_{|\goth
   p^{\lambda_\psi-1}}\not \equiv 1 .
\end{equation*}.
\subsection{} Let $X$ be an $l$-space, that is,  a Hausdorff topological   space in which  every point has a fundamental system of open compact neighborhoods. For $S$   a complex vector space, we denote by  $\mathcal C^{\infty}(X, S)$ the space  of locally constant functions $\var : X \longrightarrow S$ and by  $\mathcal C_c^{\infty}(X, S)$ the subspace consisting of those functions with  compact support.  When  $S = \cit$, we simply write $\mathcal C^{\infty}(X)$ and $\mathcal C_c^{\infty}(X)$, respectively.

 Let  $V$ be  a finite dimensional  vector space  over $F$, and 
let $\diff_Vv$ be a Haar measure on $V$. For $\var \in  \mathcal C_c^{\infty}(V)$, we define its Fourier transform by: 
$$ \mathcal F (\var) (l)= \widehat{\var} (l)= \int_V \var(v)  \psi(\langle l, v\rangle)\diff_Vv, \, \, l \in V^*,       $$
where $V^*$ denotes the space of all linear forms on $V$. It is well known that  the Fourier transform $\mathcal F$ defines  an isomorphism from $\mathcal C_c^{\infty}(V)$ onto  $\mathcal C_c^{\infty}(V^*)$. Moreover,  there exists a Haar measure $\diff_{V^*}l$ on $V^*$ such that: 
$$     \var(0)= \int_{V^*}    \mathcal F (\var) (l)  \diff_{V^*}l, \,\,   \var \in \mathcal C_c^{\infty}(V).     $$
The Haar  measure  $\diff_{V^*}l$ is called the dual Haar measure of $\diff_Vv$.
\subsection{}
Let $\bf {G}$ be a unipotent algebraic group defined over  $F$. The set of rational points $G$ of $\bf {G}$  is a unipotent $p$-adic group  with Lie algebra  $\g$ and the exponential map $\exp$ realizes   a diffeomorphism from  $\g$ onto  $G$. 

The ascending central series of $\g$ is the sequence of ideals
\begin{equation}\label{ascen}
Z^0(\g)=0  \subset  Z^1(\g)  \subset Z^2(\g) \subset \ldots   \subset Z^r(\g)= \g
\end{equation}
 defined inductively as follows   : 
for $k\geq 0$,   $Z^{k+1}(\g)$  is the inverse image in $\g$ of the center of the quotient Lie algebra 
$\g /   Z^{k}(\g)  $ by the canonical projection $\g \longrightarrow  \g /   Z^{k}(\g) $. We have 
$$    [\g, Z^{k+1}(\g)] \subset Z^{k}(\g) , \forall k.   $$
Fix a basis $(X_1, X_2, \ldots, X_n)$ of $\g$ adapted to  the ascending central series $\eqref{ascen}$ such that the lattice  
$$   \a_0 = \O-span\{X_1, X_2, \ldots, X_n\}             $$
satisfies 
$$   [\a_0,\a_0] \subset \varpi^{\beta }\a_0,       $$
where $ \beta =    [\frac{3 v(p)}{2(p -1 )}] + 1$ unless $p = 3$, when $ \beta =
v(3) + 1$ (see \cite[page 366]{howe-1977} for more details).     

For each $j \in \N\setminus\{0\}$, set 
$$  \a_j = \O-span\{\varpi^{-2^{n-1}j\beta} X_1,\varpi^{-2^{n-2}j\beta} X_2, \ldots, \varpi^{-j\beta}X_n\}     .$$
We obtain an increasing sequence of lattices in $\g$ : 
	\begin{equation}\label{latt}
\a_0 \subset \a_1 \subset \a_2\subset \ldots, \quad \g = \cup_{j \in \N}\a_j, \quad [\a_j, \a_j]\subset                \varpi^{\beta } \a_j 
\end{equation}
Put $$A_j = \exp(\a_j), \, \, j \in \N.       $$
Then,  we  have  an increasing sequence of open compact subgroups in $G$ :
\begin{equation}\label{latti}
A_0 \subset A_1 \subset A_2\subset \ldots, \quad G = \cup_{j \in \N}A_j
\end{equation}
\subsection{}\label{tang}

If $H $ is an algebraic subgroup of $G$ with Lie algebra $\h$, and $\diff_Hx$ (resp. $\diff_\h X$) denotes  a Haar measure on $H$ (resp. $\h$), we say that these measures are tangent if they satisfy 
$$    \int_H \phi(x)  \diff_Hx = \int_\h \phi \circ \exp (X) \diff_\h X, \, \forall \phi \in \mathcal C_c^\infty(H).       $$ 
\section{ Coadjoint orbits and representations :  }\label{elem}
\subsection{}  For  $g_0 \in G$, we denote by  $\varphi_{g_0}$  the inner  automorphism of   $G$ defined  by : 
 \begin{center}
$\begin{array}{ccccc}
 &\varphi_{g_0}:& G &\longrightarrow& G \\
& & g &\longmapsto& g_0gg_0^{-1}
\end{array}$
\end{center}
The differential  of  $\varphi_{g_0}$   at the identity element  induces an   automorphism of the   Lie algebra $\g$ denoted by $\mathrm{Ad}g_0$. We have 
$$     g_0 \exp(X)g_0^{-1} = \exp(\mathrm{Ad}g_0.X), \, \, \forall X \in \g.   $$
The map $\mathrm{Ad}: G \longrightarrow \mathrm{GL}(\g)$ is a  morphism of  Lie groups, called the adjoint representation of   $G$. 

 \subsection{}
 We denote by  $\g^*$ the dual  vector space of $\g$.
 The  group  $G$ acts on   $\g^*$  by the  coadjoint  representation :
$$   \langle  \mathrm{Ad}^*g.l,X  \rangle = \langle  l,\mathrm{Ad}g^{-1}.X \rangle,\, \,  g\in G, l\in \g^*, X\in \g.                      $$
  The set 
$ G.l =\left\{ \mathrm{Ad}^*g.l \, \vert \,   g\in G\right\}$
is called the  coadjoint orbit of  $l$  in $\g^*$. We denote by  $G \backslash \g^*$ the set of  all coadjoint orbits.

\noindent For $ l \in \g^*$, let    $\beta_l$ denote  the skew-symmetric bilinear form  on  $\g$  given by : 
\begin{equation}\label{bil}
 \beta_l(X,Y ) = \langle l, [X,Y] \rangle, \, \,  X,Y \in \g                 
\end{equation}
and by    $ \g(l) = \{X \in \g \, \vert \,  \langle l, [X,Y] \rangle =0, \, \forall \,  Y \in \g \}$ its radical.
 A subalgebra $\h$ of  $\g$ is said to be  subordinate  to  $l$ if  $\h$ is a totally  isotropic subspace  for 
	$\beta_l$,  that is, $ \beta_l(\h,\h) =0$. If, in addition, $\h$ is of maximal dimension,  it   is  called a   polarization at  $l$.

\subsection{}\label{reali} A unitary representation  of $G$ in a Hilbert space $\mathcal H_\pi$   means a homomorphism $\pi$ of $G$ into the group $U(\mathcal H_\pi)$ of unitary operators on $\mathcal H_\pi$ such that for every $v \in \mathcal H_\pi$ the map $G \ni x\longmapsto \pi(x)v \in \mathcal H_\pi$ is continuous. 

\noindent A subspace $\mathcal L$ of $\mathcal H_\pi$ is said to be invariant under $\pi$ if
$\pi(x) \mathcal L \subset \mathcal L$, for every $x \in G$. The representation $\pi$ is called irreducible if the only closed $\pi$-invariant subspaces of $\mathcal H_\pi$ are ${0}$ and $\mathcal H_\pi$.

\noindent Let $\rho$ be another  unitary representation  of $G$ in a Hilbert space $\mathcal H_{\rho}$. 
The representations $\pi$ and $\rho$  are said to be equivalent if there exists a unitary isomorphism
 $T : \mathcal H_\pi \longrightarrow \mathcal H_{\rho}$ such that $T \pi(x)= \rho(x)T$ for every $x \in G$.   We then  write $ \pi \simeq \rho$.

 \noindent The set of all equivalence classes of irreducible unitary representations of $G$ is denoted by $\widehat{G}$. 

\subsection{}
\noindent Take now a closed subgroup  $K$  of $G$ and  let $\sigma$ be a unitary representation of $K$ acting in a  Hilbert space $\mathcal H_\sigma$.  Let 
\begin{eqnarray*} 
  \mathcal C_c(G, K, \sigma) = \left\{ \begin{matrix} \phi: G \longrightarrow \mathcal H_\sigma  \mbox{  continuous with compact  support modulo $K$ }\, 	 \\
 \mbox{satisfying }\phi(gh) = \sigma(h^{-1}) \phi(g), \,   h\in K, \,  g \in G 
\end{matrix} \right\}     \end{eqnarray*}
and let $\mathcal H$ be its completion  
with respect to the norm 
$$ || \var      || = \left(\int_{G/K}||\var(g)||^2_{ \mathcal H_\sigma}  \diff_{G/K}\dot g \right)^{\frac{1}{2}}          $$
where $\diff_{G/K}\dot g $ is a  $G$-invariant Borel measure on the homogeneous space $G/K$.
The representation of $G$ on $\mathcal H$ by left translations is called  the induced representation  of $G$ by
 $\sigma$ and is denoted by  $\mbox{Ind}_K^G\sigma$ : 
$$   \mbox{Ind}_K^G\sigma(g) (\phi)(g') = \phi(g^{-1}g'), \,  g,g' \in  G, \, \phi  \in \mathcal H.  $$
 Let  $\b$ be a subalgebra of  $\g$ subordinate   to 
 $l\in \g^*$.  Designate by   $\chi_l$ the character of  $B = \exp(\b)$  defined by : 
\begin{equation}\label{carac}
     \chi_l(\exp(X)) = \psi_F(\langle l, X\rangle)   , \, \,  X \in \b. 
		\end{equation}
		Denote by $\pi_{l,\b}$ the  induced representation of  $G$ : 
$     \pi_{l,\b} = \mbox{Ind}_B^G \chi_l$  acting in the Hilbert space  $\H_{l,\b}$.       
\begin{tho}\label{Kir}
\
\begin{enumerate}
	\item The  representation $\pi_{l,\b}$ is irreducible if and only if  $\b$ is a polarization at~$l$.
	\item Let $\pi$ be  an  irreducible unitary representation of $G$. Then there exists $l\in \g^*$  and a  polarization $\b$ at $l$  such that $\pi$ is equivalent to   $\pi_{l,\b}$.
	\item Let $l_1$, $l_2$  be two elements  of $\g^*$ and $\b_1$, $\b_2$ be two polarizations at  $l_1$ and $l_2$ respectively.  The representations $\pi_{l_1,\b_1}$ et $\pi_{l_2,\b_2}$  are equivalent if and only if $l_1$ and  $l_2$ belong to the same  $G$-coadjoint orbit.
\end{enumerate}
\end{tho}
This theorem  is due to  A. Kirillov \cite{kirillov-1962} for the real case. 
 N.  Matringe (\cite[Theorem 3.6]{matringe-2020} and \cite[Appendix]{matringe-2024})  gave a complete proof  of the theorem for the smooth irreducible representations.

\subsection{}\label{repst}
Let $(\pi, \mathcal H)$ be a unitary representation of $G$. A vector
$v \in   \mathcal H $ is said to be smooth with respect to $\pi$ if its stabilizer $\{x \in G \, \vert \pi(x)v = v\}$ is an  open subgroup of $G$. Let $ \mathcal H^\infty$ denote  the space  of all smooth vectors of $\pi$. It is readily  seen that $ \mathcal H^\infty$ is   a $G$-invariant subspace and is dense in $\mathcal H$. Thus, we obtain a smooth representation of $G$ acting on $ \mathcal H^\infty$, which we denote  by $\pi_\infty $.

\noindent {\bf Fundamental case :} Let $l \in  \g^*$, let $\b$ be a  polarization at $l$, and set  $B = \exp(\b)$.  Let  $\chi_l$ be   the character   of  $B $  defined by formula \eqref{carac}. Consider the irreducible unitary representation $\pi= \pi_{l,\b} $  acting in the Hilbert space  $\H= \H_{l,\b}$.          

From \cite[Corollary 4.5.1]{maktouf-2026},    $  \H^{\infty} $ is the space of functions $\phi : G\longrightarrow \cit$ such that 
\begin{itemize}
	\item $\phi(xb) = \chi_l(b^{-1}) \phi(x)$, \, for $ b\in B, \,  x \in G $ ;
	\item  $\exists K$ an open compact subgroup of $G$ such that $\phi(x_0x) = \phi(x), \, \forall x \in G, x_0 \in K$ ; 
	\item $\mbox{Supp}(\phi)$  is compact modulo $B$.
\end{itemize}
\subsection{} \label{fell} 
By Theorem \ref{Kir}, for each $g \in \g^*$, we denote by $\pi_g$ the  irreducible unitary  representation of $G$ associated with $g$. It follows that the map
\begin{equation}\label{kir}
\Pi:  \g^*\longrightarrow \widehat{G}, g\longmapsto \pi_{g} 
\end{equation}
is $G$-equivariant and factors through a  bijection $ G \backslash \g^*\longrightarrow \widehat{G}$.

\noindent  We  equip $\widehat{G}$  with the Fell topology,  defined as follows: 
For an irreducible
 representation $(\pi, \mathcal H)$  of $G$,  choose $n \in \N\setminus\{0\}$, vectors $v_1, \ldots, v_n \in  \mathcal H $,   a compact set $K \subset G$,  and a real number $\e >0$, and define
 $ \mathcal U( \pi, \mathcal H, K, v_1, \ldots, v_n, \e) \subset \widehat{G}    $
to be the set of isomorphism classes $(\rho, \mathcal K) \in  \widehat{G}$ such that there exist $w_1, \ldots, w_n \in \mathcal K$
 with the property
$$   |\langle \rho(x) w_i, w_j      \rangle - \langle \pi(x) v_i, v_j      \rangle|_\cit < \e, \, x \in K, 1 \leq i,j \leq n.          $$
Sets of the form $ \mathcal U( \pi, \mathcal H, K, v_1, \ldots, v_n, \e)$ are defined to be a basis of neighborhoods of
the point $(\pi, \mathcal H)$ in $\widehat{G}$.

\begin{lem}\label{reg} Let $(\pi, \mathcal H)$ be  an irreducible unitary 
 representation  of $G$. Let $v$ be a nonzero vector in $ \mathcal H $. The family of sets  
$$  \mathcal U( \pi, \mathcal H, K, v, \e) ,          $$
where $\e>0 $  and $K$ is a compact open  subgroup in $G$, forms a  neighborhood basis at  $(\pi, \mathcal H)$ in $\widehat{G}    $.
\end{lem}
\begin{Proof}   Let $n \in \N\setminus\{0\}$,  $v_1, \ldots, v_n \in  \mathcal H $,   $B$  a compact set in $ G$ ,  and  $\e >0$. Consider the neighborhood $  \mathcal U( \pi, \mathcal H, B, v_1, \ldots, v_n, \e) $ of 
$(\pi, \mathcal H)$. We will show that it contains a neighborhood of $(\pi, \mathcal H)$   of the form described in  the statement.

 Let $Z_{\mathcal H}$ be a subset in $\mathcal H$.
\begin{itemize}
	\item Step 1 : Assume that $Z_{\mathcal H}$ is dense in $\mathcal H$.

  Choose $\e'>0$ such that 
$$   3\e'(\e' + 1 + \sup_{1\leq i\leq n}||v_i||)< \e. $$
By assumption, there exist  $v_1', \ldots, v_n' \in  Z_{\mathcal H} $
such that 
$$||v_i - v_i'|| \leq \e', \, \, i=1, \ldots, n.$$
 It follows that we have
$$ \mathcal U( \pi, \mathcal H, B, v_1', \ldots, v_n', \e')  \subset \mathcal U( \pi, \mathcal H, B, v_1, \ldots, v_n, \e)  .  $$
\item Step 2 :  Assume that the linear subspace $span(Z_{\mathcal H})$ is dense in $\mathcal H$. By step 1, we may assume that $v_1, \ldots, v_n \in span(Z_{\mathcal H})$. Then there exist $w_1, \ldots, w_m$ vectors in $Z_{\mathcal H}$ such that 
$$     v_i = \sum_{j=1}^m\lambda_{ij}w_j , \, \, 1\leq i \leq n.                $$
Set $\delta= \max_{i, j}|\lambda_{ij}|_\cit $ and $\e' = \e m^{-2}\delta^{-2}$. Then, we have 
$$       \mathcal U( \pi, \mathcal H, B, w_1, \ldots, w_m, \e')  \subset \mathcal U( \pi, \mathcal H, B, v_1, \ldots, v_n, \e)  .               $$
\item Step 3 : Since  $(\pi, \mathcal H)$ is  irreducible then the linear span of   
 $   Z_{\mathcal H} = \{\pi(g)v \, \vert \, g \in G\} $ is dense in $  \mathcal H $. By step 2, we may assume that $v_1 = \pi(g_1)v, \ldots, v_n = \pi(g_n)v$, where $g_1, \ldots, g_n$ are elements in $G$. By 
\eqref{latti}, we can find
   a compact subgroup $K$ containing $B$ and $g_1, \ldots, g_n$. Then, we have  
$$   \mathcal U( \pi, \mathcal H, K, v, \e)  \subset  \mathcal U( \pi, \mathcal H, B, v_1, \ldots, v_n, \e) .    $$
\end{itemize}
Thus the result.
\end{Proof}
\begin{lem}\label{regsm} Let $(\pi, \mathcal H)$ be  an irreducible unitary 
 representation  of $G$. Let $v$ be a  nonzero smooth vector in $ \mathcal H $, $K$  a compact open  subgroup in $G$, and $\e>0$. Let $(\rho, \mathcal K) \in  \widehat{G}$. 

If $ (\rho, \mathcal K) \in  \mathcal U( \pi, \mathcal H, K, v, \e)$, then there exists a smooth vector  $w \in  \mathcal K^\infty $ such that 
$$  |\langle \rho(x) w, w    \rangle - \langle \pi(x) v, v      \rangle|_\cit < \e,\, \, \forall x \in K.       $$
\end{lem}
\begin{Proof}
Suppose that $ (\rho, \mathcal K) \in  \mathcal U( \pi, \mathcal H, K, v, \e)$. By definition, there exists $w \in \mathcal K$ such that 
$$  |\langle \rho(x) w, w    \rangle - \langle \pi(x) v, v      \rangle|_\cit < \e,\, \, \forall x \in K.       $$
Since $v \in  \mathcal H^\infty$, there exists a compact open subgroup  $K_0$ contained in $K$  such that 
$$    \pi(x)v = v, \, \forall x \in K_0.    $$
Let $\diff_{K_0}x$ denote  the  normalized Haar measure on $K_0$ and set 
$P_{K_0}^\rho         =\int_{K_0}\rho(y)  \diff_{K_0}y$, the orthogonal projection onto 
$$\mathcal K^{K_0} : = \{ u \in\mathcal K \, \vert \, \rho(x)u=u, \forall x \in K_0  \}.$$
 Let  $ x\in K$. For all  $y, z \in K_0$, we have 
\begin{eqnarray*}
|\langle \rho(x)\rho(y) w, \rho(z) w      \rangle - \langle \pi(x) v, v   \rangle|_\cit  &=&
 |\langle \rho(z^{-1}xy) w, \rho(z) w      \rangle - \langle \pi(z^{-1}xy) v, v   \rangle|_\cit < \e. 
\end{eqnarray*}
Integrating the above inequality with respect to $\diff_{K_0}y\diff_{K_0}z$ , we obtain 
$$
| \langle \rho(x)P_{K_0}^\rho (w), P_{K_0}^\rho  (w)      \rangle  - \langle \pi(x) v, v   \rangle|_\cit  
 < \e. 
$$
Thus the vector $P_{K_0}^\rho (w)$  has the required property, and the conclusion of the lemma follows.  
\end{Proof}

\subsection{} Keep the notations of section \ref{fell}. Denote by $ \widehat{G}_s $  the smooth dual of $G$, that is, 
the set of isomorphism classes of irreducible smooth representations of $G$. We equip   $ \widehat{G}_s $ with the Fell topology defined in the same way  as for $ \widehat{G}$.
\begin{pro}\label{smooth}
The map $$   \widehat{G} \longrightarrow  \widehat{G}_s ,  (\pi, \mathcal H)  \longmapsto (\pi_{\vert}, \mathcal H^\infty)   $$
is well defined and is a homeomorphism. 
\end{pro}
\begin{Proof}
It is an immediate consequence of Lemmas \ref{reg}-\ref{regsm} and Theorem \ref{Kir}.
\end{Proof}
\begin{tho} 
 The bijection $ G \backslash \g^*\longrightarrow \widehat{G}$ is a homeomorphism (The set
of coadjoint orbits $ G \backslash \g^*$ is equipped with the quotient topology  induced by  the natural topology on $\g^*$).
\end{tho}
\begin{Proof}
It follows from  Proposition \ref{smooth} and  \cite[Theorem 3.1]{boyarchenko-sabitova-2008} that   this bijection is an open map. Moreover, by the same argument as in the real case, due to  A. Kirillov \cite[Theorem 8.2]{kirillov-1962}, the map $ G \backslash \g^*\longrightarrow \widehat{G}$ is 
 continuous. 
\end{Proof}

\section{Algebraic subgroups of codimension $1$ :}\label{codim}
\subsection{}\label{ideal}
Let $\g_0$ be an ideal of codimension $1$ in $\g$ and $G_0 = \exp(\g_0)$ the normal  algebraic subgroup  of $G$ corresponding to  $\g_0$. The following definition is due to  A. Kirillov : 
\begin{defi}
Let $l \in \g^*$.     We say that 
\begin{itemize}
	\item $l$ is of the first kind with respect to $\g_0$ if $\g(l) \subset \g_0$;
	\item  $l$ is of the second kind with respect to $\g_0$ if $\g(l) \not\subset \g_0$.
\end{itemize}
\end{defi}
Choose $X_0 \in \g \setminus \g_0$. Then we have $\g = \g_0 \oplus FX_0      $. 
\noindent Finally, let $p_{\g, \g_0} : \g^* \longrightarrow \g_0^*$ be the canonical projection. 

\begin{tho}\label{th1}
Let $l \in \g^*$, $l_0 = p_{\g, \g_0}(l)$, $\Omega= G.l$ and  $\Omega_0= G_0.l_0$.

\begin{enumerate}
	\item Suppose that $\g(l) \subset \g_0$, that is, $l$ is of the first kind with respect to $\g_0$. Then
	\begin{itemize}
		\item $l+ \g_0^\perp\subset \Omega$ ;
		\item $\operatorname{Ind}_{G_0}^{G}\pi_{l_0} \simeq \pi_l$ ; 
		\item $p_{\g, \g_0}(\Omega)=\bigcup_{s \in F} (\exp(sX_0).\Omega_0)$
	\end{itemize}
	\item Suppose that  $\g(l) \not\subset \g_0$, that is, $l$ is of the second kind with respect to $\g_0$. Then
	\begin{itemize}
		\item  $(l+ \g_0^\perp)\cap  \Omega= \{l\}$  and $p_{\g, \g_0}(\Omega) = \Omega_0$.
		\item Assume $X_0 \in \g(l)\setminus \g_0(l_0)$.  For each $t \in F$, let  $l_{(t)}$ be the element of $\g^*$ given  by 
		$$      p_{\g, \g_0}(l_{(t)})= l_0 \mbox{ and } l_{(t)}(X_0)= t.       $$
		Then we have 
		$$      \operatorname{Ind}_{G_0}^{G}\pi_{l_0} \simeq \int^\oplus_F\pi_{l_{(t)}}  \diff_Ft .    $$
	\end{itemize}
\end{enumerate}
\end{tho}
\subsection{} Let $f \in \g^*$ and $V$ be a vector subspace of $\g$. 
The following result is of particular interest.
\begin{pro}\label{lez} Suppose that $\g(g)$ is a proper vector subspace of $\g$ for at least one linear form $g \in f + V^\perp$. 	

Then there exist a nonempty Zariski open subset $U$ and a proper  Zariski closet  subset $W$  of $f + V^\perp$ such that one of  the following conditions holds: 
$$ \text{either}\, \,  \g(g) \subset \g_0, \, \, \forall g \in U$$
$$\text{ or }\, \,    \{g \in U\, \vert \, \g(g) \subset \g_0 \}  \text{ is contained in }  W  .          $$
\end{pro}
\begin{Proof}
The proof is based on a parametrization of $\g(g)$. 
Fix a basis  $(X_1, \ldots, X_n)$ of $\g$ such that $(X_1, \ldots, X_{n-1})$  is a basis of $\g_0$.
Denote by $ (X_1^*, \ldots, X_n^*) $ its dual basis.

For $g \in  f + V^\perp$, write 
$$  ad^*X_i.g    =\sum_{j=1}^n \lambda_{i,j}(g) X_j^*,      $$
where $ \lambda_{i,j}(g)= -\langle g, [X_i, X_j]  \rangle $.

Let $Y= \sum_{i=1}^n y_i X_i \in \g$. We have 
\begin{equation}\label{eq1}
  ad^*Y.g    = \sum_{j=1}^n \left(  \sum_{i=1}^n  \lambda_{i,j}(g) y_i         \right) X_j^*         
\end{equation}
It follows that 
 \begin{equation}\label{eq2}
 Y \in \g(g) \iff  \, \,   \sum_{i=1}^n  \lambda_{i,j}(g) y_i     =0, \forall j=1, \ldots, n          
\end{equation}
Moreover 
\begin{equation}\label{eq3}
 \g(g) \subset \g_0 \iff  \, \, \forall  Y \in \g(g), \, \, y_n=0.         
\end{equation}
Now denote by $M(g) =( \lambda_{i,j}(g))_{1\leq i,j\leq n }$, the  $(n\times n)$-matrix. By  assumption, 
there exists $l\in f + V^\perp$ such that $\g(l)\subsetneq \g$.
From \eqref{eq1}, we see   that $M(l)$ is   nonzero. 	Let $d \in  \N\setminus\{0\}$ be the maximum rank of $M(g)$ as $g$ ranges over  $f + V^\perp$.  Choose a $(d\times d)$-submatrix   $A(g)$ of $M(g)$  whose  determinant is  not identically zero.   Denote this determinant by $m(g)$. Thus $m(g)$,   is polynomial function on $g$. Set
$$  U= \{g \in f + V^\perp \, \vert \, m(g) \neq 0\}  .     $$
Then $U$ is a  nonempty Zariski open subset of  $f + V^\perp$. 

We now solve the linear system \eqref{eq2} over $U$ for the variables $y_i$, using the square submatrix $A(g)$ : 
Assume that the non-pivot variables, which we regard as free variables,   are $ y_{i_1}, \ldots,  y_{i_{n-d}}     $. Using Cramer's rule, we obtain
$$    y_i= \frac{P_i( y_{i_1}, \ldots,  y_{i_{n-d}}, g)}{m(g)}, \, \, i\in \{1, \ldots, n\} \setminus  \{ {i_1}, \ldots,  {i_{n-d}}\},             $$
where $P_i( y_{i_1}, \ldots,  y_{i_{n-d}}, g)$ is a polynomial function  of $y_{i_1}, \ldots,  y_{i_{n-d}}$,
with coefficients that are polynomial functions  of $g$ on $U$.

Suppose that 
$$ \{g \in U\, \vert \, \g(g) \subset \g_0 \} \subsetneq U. $$ 
Let $g \in U$. We have $ \g(g) \subset \g_0       $ is equivalent to  $y_n=0$, for every $Y\in \g(g)$. Hence,  $y_n$ cannot be a free variable and must therefore be one of the   pivot variables. Thus
$$
  y_n = \frac{P_n( y_{i_1}, \ldots,  y_{i_{n-d}}, g)}{m(g)} = 0, \forall y_{i_1}, \ldots,  y_{i_{n-d}} \in F.
$$
Since $m(g)\neq 0$ on $U$, the above condition is equivalent to the vanishing of all the coefficients $\mathcal I$ of the polynomial $P_n(y_{i_1},\ldots,y_{i_{n-d}},g)$, viewed as a polynomial in $(y_{i_1},\ldots,y_{i_{n-d}})$. Hence 
$$   \{g \in U\, \vert \, \g(g) \subset \g_0 \} =
 \{g \in U \, \vert \, P(g)= 0, \, \forall P \in {\mathcal I}\}.           $$

The result follows.
\end{Proof} 
\begin{pro}\label{negl}
 Let $\mu_{f,V}$ denote the translate by $f$ of a Haar measure on $V^\perp$.

Then 
\begin{itemize}
	\item Either the set of linear forms in $      f + V^\perp $ of the first kind with respect to $\g_0$ is $\mu_{f,V}$-null.
	\item  Or the set of linear forms in $      f + V^\perp $ of the second kind with respect to $\g_0$ is $\mu_{f,V}$-null.
\end{itemize}
\end{pro}
\begin{Proof} The result follows  immediately  from  Proposition \ref{lez} and the  fact that $\mu_{f,V}$ is a regular Borel measure on  $      f + V^\perp $.
\end{Proof}

\section{ The affine space $f+ \h^\perp$  :  }{} 
\

Throughout what follows,  $H$ is a proper  algebraic subgroup  of $G$ with  Lie algebra $\h$,  and   $\chi$  is  a unitary character of $H$.

Let $f$ be  a linear form on $\g$ 
such that $\chi=\chi_{f}$, where
\begin{equation*}
\chi_{f}(\exp X)=\psi_F(\langle f,X\rangle)\mbox{, }X\in\mathfrak{h}.
\end{equation*}
The Lie algebra $\mathfrak{h}$ must be subordinate to $f$, that is,
$f([\mathfrak{h},\mathfrak{h}])=0$.
Consider  the affine subspace  $X_\chi = f + \h^\perp$ of $\mathfrak g^*$ consisting of linear forms whose restriction to  $\h$ coincide with that  of $f$.

 Denote by  $Q(f, \g)$   the set of subalgebras $\goth k$  which  is subordinate to $f$  and that  ${\goth k} + \g(f)$ is a maximal   isotropic subspace (for $\beta_f$). 
Finally, denote by 
\begin{equation}
 \Xi= \{l \in f+ \h^\perp \, \vert \, \h \in Q(l, \g) \}\}
\end{equation}

\begin{lem} \label{dim} Let $ l\in  f+ \h^\perp $. We ave 
$$   \dim G.l - \dim H.l =  \dim(\g.l \cap \h^\perp) \geq   \dim H.l   .   $$
Moreover,   equality  holds if and only if $ l \in    \Xi  $.
\end{lem}
\begin{Proof} Since $  \g.l = \g(l)^\perp$, then we have  
\begin{eqnarray*}
 \dim(\g.l \cap \h^\perp)&=&  \dim(\g(l)^\perp \cap \h^\perp)= \dim(\g(l) + \h)^\perp = \dim\g - \dim(\g(l) + \h) \,\, \,  (*)\\
&=&  (\dim\g - \dim\g(l)) - (\dim \h- \dim \h(l))\\
&=& \dim G.l - \dim H.l.
\end{eqnarray*}
As $\g(l) + \h$ is isotropic with respect to  $\beta_l$  then 
$$  \dim(\g(l) + \h)\leq \frac{1}{2}(\dim \g + \dim \g(l)).    $$
Using $(*)$, we deduce that : 
$$
 \dim(\g.l \cap \h^\perp)\geq   \dim\g - \frac{1}{2}(\dim \g + \dim \g(l))= \frac{1}{2}(\dim \g - \dim \g(l)) = \frac{1}{2} \dim G.l.
$$
Moreover, one has 
\begin{eqnarray*}
\dim(\g.l \cap \h^\perp)= \dim H.l &\iff & \dim(\g.l \cap \h^\perp)=  \frac{1}{2} \dim G.l \\
&\iff& \dim(\g(l) + \h)=\frac{1}{2}(\dim \g + \dim \g(l)) \iff  l \in    \Xi.
\end{eqnarray*}
Thus the result.
\end{Proof}
Fix an ideal  $\g_0$  of codimension $1$ in $\g$ containing $\h$. So we use the notations of section \ref{ideal}.
For $l \in \g^*$, we set $l_0 = p_{\g, \g_0}(l)$.  In the sequel, the notation $\perp_0$ indicates that the orthogonal is taken in $\g_0^*$.
\begin{pro}\label{int}
 Let $ l\in  f+ \h^\perp $.  Assume that  $G.l \cap (f + \h^\perp)$  is a smooth variety of dimension equal to     $\dim(\g.l \cap \h^\perp)$  and that $ \dim(\g.l \cap \h^\perp)= \dim(\g.l' \cap \h^\perp)$, for every $l'  \in G.l \cap (f + \h^\perp)$.

Then $G.l_0 \cap (f_0 + \h^{\perp_0})$ is a  smooth variety of dimension equal to     $\dim(\g.l_0 \cap \h^{\perp_0})$.
\end{pro}
\begin{Proof} 
$\bullet$ {\it First case} : Assume in this case that $\g(l) \subset \g_0$. Then
$$ \g(l) \subset \g_0(l_0)  \mbox{ and } \dim   \g(l)  = \dim \g_0(l_0) -1    .$$ 
Fix $X_{m+1} \in  \g_0(l_0)\setminus  \g(l) $. 
 Let    $(X_1, \ldots, X_m) $ be  an  adapted complementary  basis  of $\g_0(l_0)$ in $\g$. 
Then  $(X_1, \ldots, X_m, X_{m+1}) $ is  an  adapted complementary  basis  of $\g(l)$ in $\g$.  

For $t= (t_1, \ldots, t_m, t_{m+1}) \in F^{m+1}$, set 
$$g_t= \exp(t_1X_1)\times \cdots \times \exp(t_mX_m) \exp( t_{m+1}X_{m+1}).         $$
Since the map $$ \var :   F^{m+1}    \longrightarrow G.l, t \longmapsto g_t.l        $$
is a diffeomorphism and  $G.l \cap (f + \h^\perp)     $  is a smooth subvariety, its inverse image 
$$ S= \{(t_1, \ldots, t_m, t_{m+1}) \in F^{m+1} ; \, \,  g_t.l \in f + \h^\perp\}          $$
 is a  smooth subvariety of $F^{m+1}$ of dimension equal to  $\dim(\g.l \cap \h^\perp)$. 

Let denote by $S_0 = \{s=(t_1, \ldots, t_m) \in F^{m} ; \, \,   (s, 0) \in S\} $. Then, we have 
$$ S = S_0 \times F.       $$
It is straightforward to verify that the smooth map 
$$  \psi : F^m\longrightarrow    F^{m+1} , s\longmapsto (s,0),        $$
 is transverse to $S$. Hence, by the transverse preimage theorem \cite[Theorem 6.30]{lee-2013},  $\psi^{-1}(S)$ is a smooth subvariety of $F^m$ whose codimension is equal to the codimension of $S$ in $F^{m+1}$.

Consider now the smooth map
$$    \var_0: S_0 \longrightarrow G.l_0, \,  s\longmapsto g_{(s,0)}.l_0   .   $$
We have  $ \var_0 (S_0)= G.l_0 \cap (f_0 + \h^{\perp_0}) $.
We claim that $\var_0$ 
is a smooth embedding. Indeed, using the fact that $G(l_0)= G_0(l_0)$, we see that $\var_0$ is an immersion injective.  Let $\mathcal W$ be any open subset of  $ S_0 $. Then 
  $ \mathcal W \times F     $ is open in $S$. It follows that 
	$$\var_0(\mathcal W )= p_{\g, \g_0} \circ\var(\mathcal W \times F  ) $$ is open in $G.l_0 \cap (f_0 + \h^{\perp_0})$. This proves the claim. 

 Consequently  $G.l_0 \cap (f_0 + \h^{\perp_0})= \var_0 (S_0) $ is a smooth subvariety of  $f_0 + \h^{\perp_0}$.
\begin{remm}
To show that $G(l_0)= G_0(l_0)$ we proceed as follows. Observe that the linear map 
$$ F^m\longrightarrow \g.l_0, \, (t_1, \ldots, t_m) \longmapsto \left(\sum_{k=1}^mt_kX_k\right).l_0        $$
is surjective. Hence 
$$  m-1 = \dim(G_0.l_0)  <   \dim(G.l_0) \leq m.         $$
Consequently  $  \dim(G.l_0) = m $ and that $ \dim\g(l_0)= \dim\g_0(l_0)$. Thus $ \g(l_0)= \g_0(l_0)   $
and  $G(l_0)= G_0(l_0)$.
\end{remm}
$\bullet$ {\it Second case } : Assume in this case that $\g(l) \not\subset \g_0$. Then
$$ \g_0(l_0) \subset \g(l)  \mbox{ and } \dim   \g(l)  = \dim \g_0(l_0) +1    .$$ 
The map : 
$$ \theta :  G.l \longrightarrow G_0.l_0 , l'\longmapsto  l'_{\vert \g_0} $$ 
is a diffeomorphism.  By assumption  $G.l \cap (f + \h^\perp)$  is a smooth subvariety. Therefore its image 
$$    \theta(G.l \cap (f + \h^\perp))= G.l_0 \cap (f_0 + \h^{\perp_0})$$ 
is also  a smooth subvariety of $f_0 + \h^{\perp_0}$.
\end{Proof}
\subsection{}\label{Zar}

\begin{pro}\label{Zariski}
There exists a nonempty Zariski open $\mathcal U$ of $f + \h^\perp$  such that :
\begin{enumerate}
\item $\mathcal U$ is a union of intersections with $G$-orbits;
	\item For any $l \in \mathcal U$, 
	\begin{itemize}
		\item $G.l \cap (f + \h^\perp)     $  is a smooth variety of dimension equal to     $\dim(\g.l \cap \h^\perp)$. In particular, we have $ \dim(\g.l \cap \h^\perp)= \dim(\g.l' \cap \h^\perp)        $, for all $l'  \in G.l \cap (f + \h^\perp)$.
		\item $\dim(\g.l \cap \h^\perp) = \dim H.l$  if and only if $ G.l \cap (f + \h^\perp)  \subset  \Xi$.
	\end{itemize}
	\item $\Xi \cap \mathcal U$ is  a Zariski open subset  in $f + \h^\perp$. 
\end{enumerate}

\end{pro}
\begin{Proof}
Using the notation introduced in the  Appendix, let $j$ be the smallest integer $1 \leq j \leq m$ such that 
$$  F_{j-1}\supset  f + \h^\perp \mbox{ and }    F_{j} \not\supset  f + \h^\perp  .      $$
Then 
$$  \mathcal V    :=  ( f + \h^\perp ) \cap   ( F_{j-1} \setminus F_j)= \{l \in f + \h^\perp ; \, Q_j(l) \neq 0 \}   . $$
Since $F_{j-1} \setminus F_j$ is $G$-invariant,  $\mathcal V$ is a union of intersections with $G$-orbits.

Let $ \lambda  :F_{j-1} \setminus F_j
\longrightarrow \g^* $   be  the regular rational function  :
$$ \lambda(x)= \sum_{k=1}^m  P_k(0,x) v_k   .           $$
(See appendix, Corolloary \ref{AC}) and
let denote by $ \lambda_0  :
  \mathcal V  \longrightarrow \g^* $  its restriction to $\mathcal V $. Then $\lambda_0$ is a  regular rational function on $\mathcal V$.
   Denote by $  \mathcal U $ the Zariski-open subset of  $ \mathcal V $ on which  the differential of  $\lambda_0$ has maximal  rank. Then  $  \mathcal U $ is a union of intersections with $G$-orbits in $ f + \h^\perp$.

Now the restriction $\lambda_{0\vert \mathcal U} $ is a subimmersion. 
We deduce that, for any $l \in \mathcal U$, $ \lambda_0^{-1}(\lambda_0(l))$ is a smooth subvariety in
 $  f + \h^\perp      $ and  $T_{l} \lambda_0^{-1}(\lambda_0(l))$, the tangent subspace of $\lambda_0^{-1}(\lambda_0(l))$ at $l$, is given by 
$$       T_{l} \lambda_0^{-1}(\lambda_0(l)) = \ker d_l \lambda_0 = \ker d_l \lambda \cap \h^\perp.    $$
But, we have $  \lambda_0^{-1}(\lambda_0(l)) = G.l \cap   ( f + \h^\perp )$. Thus $  G.l \cap   ( f + \h^\perp )     $ is a  a smooth subvariety in
 $  f + \h^\perp      $. We see immediately tha $  T_l G.l \cap   ( f + \h^\perp ) \subset \g.l \cap  \h^\perp            $. 

On the other hand,  the restriction  of $ \lambda$ to $ G.l$ is a constant function. Consequently
$$\g.l = T_l G.l \subset \ker d_l \lambda    .     $$
So that $ \g.l \cap \h^\perp \subset   T_{l} \lambda_0^{-1}(\lambda_0(l))        $. Hence $   T_l G.l \cap   ( f + \h^\perp )=  \g.l \cap \h^\perp $.

Now we prove  that $\Xi\cap\mathcal U$ is a Zariski open in  $f+\h^\perp$. Let $2d$ denote the common dimension of $G.l$,  for $l \in\mathcal U$.
We have 
\begin{eqnarray*}
 \dim(\h + \g(l))&=& \dim \h +\dim \g(l) - \dim(\h\cap \g(l))\\
&\leq& \frac{1}{2}(\dim \g + \dim \g(l))
\end{eqnarray*} 
So that 
$$    \dim(\h\cap \g(l)) \geq \dim\h - \frac{1}{2}\dim G.l  =   \dim\h - d.   $$
It follows that 
 \begin{eqnarray*}
 \h \in Q(l, \g) &\iff&   \dim(\h\cap \g(l)) = \dim\h -d.
\end{eqnarray*}
But $\dim(\h\cap \g(l)) >(\dim\h - d)$ means that $\g(l) \cap \h\cap W \neq 0$ for every subspace $W$ of $\g$ of dimension $\dim \g -(\dim\h - d)$.  
From \cite[Lemma I.II.4]{dixmier-1974},  we deduce the desired assertion.  
\end{Proof}

\begin{co}\label{zar1}
 Assume that  $\Xi \cap \mathcal U \neq \varnothing$. Let $\mathcal U_0 = p_{\g, \g_0}(\Xi \cap \mathcal U )$.

Then 
\begin{itemize}
	\item $\mathcal U_0 $   is  a Zariski open subset  in $f_0 + \h^{\perp_0}$. 
	\item  For all $ l_0 \in\mathcal U _0 $,  $G.l_0 \cap (f_0 + \h^{\perp_0}) $  and   $G_0.l_0 \cap (f_0 + \h^{\perp_0}) $ are   smooth varieties of dimension equal to $\dim H.l_0$. 
\end{itemize}  
\end{co}

\begin{Proof}
The first assertion  follows from  Proposition \eqref{Zariski} and   the fact that the restriction of 
$p_{\g, \g_0}$ to $f + \h^{\perp}$ is a Zariski-open map onto $f_0 + \h^{\perp_0}$. 

Let  $ l_0 \in\mathcal U _0 $ and let $l \in \Xi \cap \mathcal U $ such that 
$    p_{\g, \g_0}(   l)= l_0$. We will compare the $H$-orbits in  $G.l \cap (f_0 + \h^{\perp}) $ and the $H$-orbits in  $G.l_0 \cap (f_0 + \h^{\perp_0}) $.

By Proposition \ref{Zariski} and \cite[Proposition 3.2.1]{maktouf-2026}, we have 
$$ G.{l} \cap (f+ \h^{\perp}) = H.{l}_1 \sqcup \ldots \sqcup  H. {l}_r .  $$
 Since $  p_{\g, \g_0}$ is $H$-equivariant, we obtain that 
$$   p_{\g, \g_0}(G.{l} \cap (f+ \h^{\perp})) = H. l_{1,0} \cup \ldots \cup  H.l_{r,0}, \, \,     $$
 where  $l_{i,0}=  p_{\g, \g_0}({l}_i )$, in particular $ l_{1,0}=l_0       $.    
We distinguish two cases : 

$\bullet$ {First case} : $\g({l}) \subset \g_0$. In this case, we have
\begin{enumerate}
	\item 
 $ \g({l}) \subset \g_0(l_0)$  and 
$ \dim   \g({l}) + 1 = \dim \g_0(l_0)     $.   
\item For all $l_0' \in G.l_0\cap (f_0+ \h^{\perp_0}) $, we have 
$$ \dim H.l_0' = \dim H.{l} -1.         $$
\item We have  $$   p_{\g, \g_0}(G.{l} \cap (f+ \h^{\perp})) = G.l_0\cap (f_0+ \h^{\perp_0})   =  H. l_{1,0} \sqcup \ldots \sqcup  H.l_{r,0} ,  $$
which is a smooth subvariety with dimension equal to $\frac{1}{2}\dim G_0.l_0$.
\end{enumerate}
Indeed, it is readily seen that
$$     \dim(\h \cap\g_0(l_0)) \leq \dim(\h \cap\g({l})) +1.        $$
It follows that :
\begin{eqnarray*} 
 \dim(\h + \g_0(l_0))&  = &\dim \h + \dim \g_0(l_0) -  \dim(\h \cap\g_0(l_0)) \\
& \geq&\dim \h + \dim \g({l})- \dim(\h \cap\g({l})) \\
&=&\dim (\h + \g({l})) = \frac{1}{2}(\dim \g + \dim \g({l}))\\
&=& \frac{1}{2}(\dim \g_0 + \dim \g_0(l_0)).
\end{eqnarray*}
But $\h + \g_0(l_0)$ is an isotropic subspace with respect to the skew-symmetric bilinear form  $\beta_{l_0}$  on  $\g_0$   defined by :
$$ \beta_{l_0}(X,Y ) = \langle {l_0}, [X,Y] \rangle, \, \,  X,Y \in \g_0.                $$
We deduce that 
$$   \dim(\h + \g_0(l_0))=  \frac{1}{2}(\dim \g_0 + \dim \g_0(l_0)).       $$
Thus $\h \in  Q(l_0, \g_0)$. Hence 
$$    \dim H.l_0 = \frac{1}{2}\dim G_0.l_0  =  \dim(\g_0.l_0\cap\h^{\perp_0})=   \frac{1}{2}\dim G.{l}-1 =  \dim H.{l}-1.      $$ 
We now show that the $H$-orbits are pairwise disjoint : let $i, j$ be such that $ H. l_{i,0} = H. l_{j,0}$. Then there exist $h \in H$ and $g \in G$ such that 
$$   h  . l_{i,0}  = g. l_{i,0} = l_{j,0} .  $$
It follows that $ h^{-1}g \in G(l_{i,0}) = G_0(l_{i,0})     $. 
We claim that $ G_0(l_{i,0}).l_i= H (l_{i,0}).l_i$.  It is immediate that  
$$H (l_{i,0}).l_i \subset G_0(l_{i,0}).l_i= l_i + \g_0^\perp.  $$      
Since $l_i + \g_0^\perp$ is a one-dimensional irreducible closed subset of $\g^*$, the inclusion is proper if and only if $H (l_{i,0}).l_i = l_i$.  However,  this would  imply $  H (l_{i,0})= H (l_{i})$ and hence $\dim H .l_i = \dim H.l_{i,0}$  which is impossible. Therefore, there exists $ h' \in H (l_{i,0})$ such that 
$$    h'. l_{i} =  h^{-1}g .l_{i}  . $$
Consequently $H. l_{i} = H. l_{j}$.

$\bullet$ {Second case} : $\g({l})\not \subset \g_0$.   In this case, we have  
\begin{enumerate}
	\item  $ \g_0(l_0) \subset \g({l}) $  and  $ \dim   \g({l}) = \dim \g_0(l_0) +1 $ ,
	\item  
	$  p_{\g, \g_0}(G.{l} \cap (f+ \h^{\perp}))   = G_0.l_0 \cap (f_0 + \h^{\perp_0})=  H. l_{1,0} \sqcup \ldots \sqcup  H.l_{1,0},$
which is a smooth subvariety  of  dimension equal to $\frac{1}{2}\dim G_0.l_0$.
\end{enumerate}
Thus the result.
\end{Proof}

\begin{co}\label{zar2} We have 
\begin{itemize}
	\item \label{a} if  $\Xi \cap \mathcal U \neq \varnothing$ then $$\dim G.l = 2 \dim H.l, \, \, \forall l \in \Xi \cap \mathcal U  ; $$ 
	\item \label{b} if $\Xi \cap \mathcal U  =  \varnothing$ then $$\dim G.l > 2 \dim H.l, \, \, \forall l \in\mathcal U .$$ 
\end{itemize}  
\end{co}

\section{Orbital  Decomposition of  $\operatorname{Ind}_{H}^{G}\chi$ : }{}\label{cent}
\

In this section, we give   the  central  decomposition of $\operatorname{Ind}_{H}^{G}\chi$.
\subsection{The decomposition of measures :}\label{meas22}
\

In the following, $K$ designates either $G$ or $H$. 
The map
$$
K\times X_\chi\longrightarrow\mathfrak g^*,
\qquad (x, l)\longmapsto x .l,
$$
is a morphism of algebraic varieties. Hence, by Chevalley's  theorem \cite[Th\'eor\`eme 1.8.4, Chapitre IV]{grothendieck-1964}, its image $C : = K.X_\chi$ is constructible in 
$\mathfrak g^*$,  for the Zariski topology.

The set $ K.X_\chi$ is locally compact, and has a countable base for   the $p$-adic topology.
Furtheremore, any $K$-orbits is closed. Then $ K.X_\chi/K$ is a $T_1$-space and  countably separated. 

The map  $\begin{array}{ccccc}
 && X_\chi &\longrightarrow& G.X_\chi/G \\
& & l &\longmapsto& G.l
\end{array}$ is continuous,  $H$-equivariant and factors to a surjective map :
$$P_\chi : X_\chi/H \longrightarrow  G.X_\chi/G.$$
Thus, the following diagram is commutative : 
\[
\begin{array}{ccc}
X_\chi &  \xrightarrow{\hspace{0.3cm}  \hspace{0.5cm}} &  X_\chi/H \\
  &  \searrow & \big \downarrow P_\chi \\
 &  &  G.X_\chi/G
\end{array}.
\]
We see immediately that, for any $\Omega\in  G.X_\chi/G$, 
$$ P_\chi^{-1}(\{\Omega\}) = \{\omega \in X_\chi/H \, \, ;\, \, \omega \subset \Omega  \cap X_\chi \}.       $$
Let $ \mu_{\h} $ be a Haar measure on $  \h^\perp$,  and let $ \mu_{f, \h}$ be its translate 
  by $f$ to $X_\chi$. 
We regard  $ \mu_{f, \h}$  naturally as a Borel measure on $\g^*$.
We  denote  by $ \nu_{f, \h}$ a  finite Borel measure on $X_\chi$ that is equivalent to 
$ \mu_{f, \h}$.

Let  $\dot\nu_{f, \h}^G$ (resp. $\dot\nu_{f, \h}^H$) be the push-forward of $\nu_{f, \h}$ under the map  $X_\chi\longrightarrow G.X_\chi/G$ (resp. $X_\chi\longrightarrow X_\chi/H$). Then $\dot\nu_{f, \h}^G$ is the  push-forward of $\dot\nu_{f, \h}^H$ under  $ P_\chi$.
Let 
 \begin{equation}\label{form20}
\dot\nu_{f, \h}^H = \int_{ G.X_\chi/G}{\mu_y}  \diff\dot\nu_{f, \h}^G(y)
\end{equation}
be  the disintegration of $\dot\nu_{f, \h}^H$ with respect to $\dot\nu_{f, \h}^G$ and the Borel map  $P_\chi$ (see \cite[\S 3]{bourbaki-2007}).
\begin{lem}\label{desint11} The notations are as above.
The support of $\mu_y$ is $  P_\chi^{-1}(y)$, $\dot\nu_{f, \h}^G$-almost every $y$.
\end{lem}
\begin{dem}
 The result follows by the same argument as that used in \cite[Lemme 1]{fujiwara1}.
\end{dem}
\begin{pro}\label{desin221} We  have 
\begin{equation}\label{form22}                  
 \int^\oplus_{X_\chi/H} \pi_l  \diff\dot\nu_{f, \h}^H(l) = \int^\oplus_{G.X_\chi/G}n_\Omega \pi_\Omega  \diff\dot\nu_{f, \h}^G(\Omega)
\end{equation}
where  $ n_\Omega= card\{\mbox{$H$-orbits in $\Omega \cap X_\chi$}\}$.
\end{pro}
\begin{dem}
We distinguish two cases:
\begin{enumerate}
	\item[(i)] Case I:  $(\Xi \cap \mathcal U \neq \varnothing)$.  
	Set $Y = P_\chi((\Xi \cap \mathcal U)/H)$. Since  $\Xi \cap \mathcal U$ is a union of intersections with
	$G$-orbit, one has $P_\chi^{-1}(G.X_\chi/G \setminus Y) =  X_\chi/H\setminus(\Xi \cap \mathcal U)/H$ which
	is a set with measure zero. Hence, we may replace  $      X_\chi/H   $  with $(\Xi \cap \mathcal U)/H$ and $ G.X_\chi/G $ with $Y$. For every $ \Omega \in Y$, the fiber  $ P_\chi^{-1}(\{\Omega\}) $ is  finite and has  cardinality $n_\Omega$. Thus,
	$$  \int_{ P_\chi^{-1}(\{\Omega\})} \pi_\omega  \diff {\mu_\Omega}(\omega) \simeq n_\Omega \pi_\Omega.   $$ 
	The result then follows from \eqref{form20}  and \cite[Theorem 2.11]{mackey-1953}.
	\item[(ii)] Case II:   $(\Xi \cap \mathcal U = \varnothing)$. Then  $\dim G.l > 2 \dim H.l, \, \, \forall l \in\mathcal U$. Consequently, $ P_\chi^{-1}(P_\chi(H.l))$ is infinite, for every $l \in \mathcal U$ and 
	$$    \int^\oplus_{P_\chi^{-1}(P_\chi(H.l))} \pi_{H.l'}  \diff\mu_{P_\chi(H.l)} \simeq \infty \pi_{P_\chi(H.l)} . $$
\end{enumerate}
Thus the result.
\end{dem}
\subsection{} The first main result of this paper is the following:
\begin{tho} \label{des} Keeping the above notation, we have
\begin{equation}\label{prin} 
   \operatorname{Ind}_{H}^{G}\chi = \int^\oplus_{G.X_\chi/G}n_\Omega \pi_\Omega  \diff\dot\nu_{f, \h}^G(\Omega)
	\end{equation}
where  $n_\Omega \in \N \cup \{+\infty\}$ is given as follows : 
\begin{enumerate}
	\item  if  $\Xi \cap \mathcal U \neq \varnothing$, then $$ n_\Omega= card\{\mbox{$H$-orbits in $\Omega \cap X_\chi$}\}, \, \, \forall \Omega \cap \mathcal U \neq \varnothing ; $$ 
	\item  if $\Xi \cap \mathcal U  =  \varnothing$, then $$ n_\Omega= +\infty, \forall\Omega.$$ 
\end{enumerate}
\end{tho}
In the real case, 
these results are originally established by Corwin, Greenleaf, and Grelaud \cite{corwin-greenleaf-grelaud-1987} and subsequently by Lipsman \cite{lipsman-1989}. Their work provides a fundamental framework for the study of the decomposition of induced representations and the associated multiplicity problem. In the present paper, we extend these results to the setting of $p$-adic unipotent groups. Along these lines, the general strategy is inspired by the real case.

The proof of Theorem \ref{des} is given in the following  subsections.
\subsubsection{} To prove Theorem \ref{des}, it suffices, by Proposition \ref{desin221}, to establish that
\begin{equation}\label{desin225}
   \operatorname{Ind}_{H}^{G}\chi =  \int^\oplus_{X_\chi/H} \pi_l  \diff\dot\nu_{f, \h}^H(l) 
	\end{equation}
		The proof of \eqref{desin225} goes   by induction on $\dim G/H$.			

Assume first  $\dim G/H =1$. This implies that $ H $ is a normal subgroup in $G$.  Thus our result is exactly the assertions of Theorem \ref{th1}.

Henceforth, suppose that $\dim G/H >1$. Let $G_0$ be a  normal  algebraic subgroup  of $G$ of codimension one and containing $H$.  Denote  its  Lie algebra by  $\g_0$, which is an ideal of $\g$. 

Using  induction in stages and  the induction hypothesis, we obtain
\begin{eqnarray*}
\operatorname{Ind}_{H}^{G}\chi &= & \operatorname{Ind}_{G_0}^{G} \operatorname{Ind}_{H}^{G_0}\chi=  \operatorname{Ind}_{{G_0}}^{G} \int^\oplus_{X^0_\chi/H} \pi_{l_0}  \diff\dot\nu_{f_0, \h}^H(l_0) =    \int^\oplus_{X^0_\chi/H}\operatorname{Ind}_{{G_0}}^{G} \pi_{l_0}  \diff\dot\nu_{f_0, \h}^H(l_0), 
\end{eqnarray*}
where $ X^0_\chi= f_0 + \h^{\perp_0}$.
          
We will now prove the following assertion : 
\begin{equation}\label{fonda1}
 \int^\oplus_{X^0_\chi/H}\operatorname{Ind}_{{G_0}}^{G} \pi_{l_0}  \diff\dot\nu_{f_0, \h}^H(H.l_0) = 
  \int^\oplus_{X_\chi/H} \pi_l  \diff\dot\nu_{f, \h}^H(H.l). 
\end{equation}
To this end, we distinguish two cases.
\subsection{Case I: Suppose that  $\Xi \cap \mathcal U \neq \varnothing$. }\label{ahm1}
\subsubsection{} In this paragraph, assume that  $\g(l) \subset \g_0$,  $\mu_{f,\h}$-almost every
 $l \in X_\chi$. 

Let $l \in \Xi \cap \mathcal U $ such that $\g(l) \subset \g_0$.  Then, by Corollary \ref{zar1}, 
$ \h \in   Q(l, \g) \cap Q(l_0, \g_0)$. We deduce that $$\dim (\h \cap \g_0(l_0))=  \dim (\h \cap \g(l))+1.   $$
Thus, for every $l' \in H$, we have $ l'+ \g_0^\perp \subset H.l$. Consequently,  
$p_{\g,\g_0}^{-1}(H.l_0)= H.l$. 

Hence, the map
$$  \Theta :    X_\chi/H\longrightarrow   X^0_\chi/H, \, \, \,  H.l \longmapsto H.l_0       $$
 induces an isomorphism from $(\Xi \cap \mathcal U)/H $ onto $\mathcal U_0/H$. Moreover,  we have the commutative diagram 
\[
\begin{array}{ccc}
f + \h^\perp &  \xrightarrow{\hspace{0.5cm} p_{\g,\g_0} \hspace{0.5cm}} & f_0+ \h^{\perp_0} \\
\alpha\, \big\downarrow & & \big \downarrow \alpha_0 \\
(f + \h^\perp)/H & \xrightarrow{\hspace{0.5cm} \Theta \hspace{0.5cm}} & (f_0+ \h^{\perp_0})/H
\end{array},
\]
where $\alpha(l)= H.l$,  $l \in f + \h^\perp $ and $\alpha_0(l_0)= H.l_0$, $l_0 \in f_0 + \h^{\perp_0}$.
 
We see that $\dot\nu_{f_0, \h}^H$ coincides with the push-forward of  $\dot\nu_{f, \h}^H$ under $\Theta$.
Hence the equality \eqref{fonda1} follows.

\subsubsection{}  Assume now  that  $\g(l) \not\subset \g_0$,  $\mu_{f,\h}$-almost every
 $l \in X_\chi$. 

By Proposition \ref{lez}, there exists a nonempty  Zariski-open subset $ \mathcal V \subset f+ \h^\perp$ such that  
$$   \g(l) \not\subset \g_0, \, \forall l \in  \mathcal V .        $$ 
We may suppose that $\mathcal V $ is $H$-invariant and $ \mathcal V + \g_0^\perp= \mathcal V $.

Fix $X_0 \in \g \setminus \g_0$. For each $l \in   \mathcal V$, there exists a unique vector 
$X_l \in \g_0$ such that 
$$ (*) \quad   (X_0+ X_l).l= 0 \, (\mbox{which is equivalent to } X_l.l = -X_0.l).          $$
We claim that  $X_l$ varies rationally  on the coordinates of    $l$ in  $ \mathcal V$.  Indeed let $X_1, \ldots, X_{n-1}$ be a basis in 
$\g_0$ and let  $e_1^*, \ldots, e_{n}^*$  be a basis in  $\g^*$. Writing
$$ X_l= \sum_{i= 1}^{n-1} x_iX_i, \, \,  l= \sum_{j=1}^n y_j e_j^*, \, X_i. e_j^* =
 \sum_{k=1}^n a_{ij}^k e_k^*,        $$  
then $(*)$ is equivalent to 
\begin{equation}\label{linsys}
   \sum_{i=1}^{n-1}\left(\sum_{j=1}^n y_j a_{ij}^k \right)x_i= -\sum_{j=1}^n y_j a_{0j}^k , \, \, 
 1 \leq k \leq n  
\end{equation}
Since the linear  system \eqref{linsys} has a unique solution, the matrix $\left(\sum_{j=1}^n y_j a_{ij}^k\right)_{ik}$ is of maximal rank, namely  $(n-1)$.

We deduce that there  exists $1 \leq i_0 \leq n$ such that the linear system \eqref{linsys} is equivalent to 
\begin{equation}\label{linsys2}
   \sum_{i=1}^{n-1}\left(\sum_{j=1}^n y_j a_{ij}^k \right)x_i= -\sum_{j=1}^n y_j a_{0j}^k , \, \, 
 1 \leq k \leq n, k \neq i_0.  
\end{equation}
Hence, by Cramer's rule, $X_l$ depends rationally on the coordinates of $l$.
\begin{remm} : Since $[\g, \g] \subset \g_0$, we have 
$$ X_l = X_{l'}, \forall l' \in p_{\g, \g_0}^{-1}(l_{\vert \g_0}). $$ 
Consequently,  when necessary, we shall also write $ X_{l\vert \g_0}$ for $X_l$. 
\end{remm}
  The rational map 
$$k : \mathcal V \longrightarrow  \mathcal V_0\times F, \quad l \longmapsto \left(l_{\vert \g_0}, l(X_0+ X_l) \right)      $$
is well defined and bijective. 
Here $  \mathcal V_0  = p_{\g, \g_0} ( \mathcal V)$.
  
	We show that  the push-forward of $\nu_{f, \h}$ under the map $k$ is equivalent to the product measure $\nu_{f_0, \h} \times \diff_Ft$. Indeed, let $\phi$ be a positive measurable function on $X_\chi$.
	Then
	$$ \int_{\mathcal V_0\times F} \phi(k^{-1}(l_0, t(X_0+ X_{l_0})) \diff\nu_{f_0, \h}(l_0) \times \diff_Ft =0     $$
	if and only if, $  \nu_{f_0, \h}$-almost every  $l_0 \in \mathcal V_0$, 
	$$   \int_{ F} \phi(k^{-1}(l_0, t(X_0+ X_{l_0}) )\diff_Ft =0,          $$
	which is equivalent to 
	$$    \phi(l)=0, \forall l \in p_{\g, \g_0}^{-1}(l_0).         $$
	By Fubini's theorem, the above condition is equivalent to
	$$    \int_{\mathcal V} \phi(l) \diff\nu_{f, \h}(l) =0.     $$
	\subsubsection{} Keep notations and assumptions  of the above section. 
	\begin{lem}\label{magic}
	For every $x \in H$ and $l \in \mathcal V$, we have
	$$ X_{x.l}= x.X_0- X_0 + x.X_l.         $$
	\end{lem}
	\begin{Proof}
	It is immediate that $  x.X_0- X_0 + x.X_l\in \g_0$, since $ x.X_0- X_0  \in  \g_0 $.
	Moreover 
	$$    \langle (x.X_0+ x.X_{l}).(x.l), Y \rangle =  -\langle  l, [X_0+ X_l,x^{-1}.Y ] \rangle=0, \forall Y \in \g.    $$
	Thus the result.
	\end{Proof}
	By Lemma \ref{magic}, the $p$-adic number $  l(X_0+ X_l)$  depends only on the $H$-orbit of $l$.
	We deduce that the following diagram is commutative : 
	\[
\begin{array}{ccc}
\mathcal V &  \xrightarrow{\hspace{0.5cm} k  \hspace{0.5cm}} & \mathcal V_0\times F\\
\alpha\, \big\downarrow & & \big \downarrow \alpha_0\times id_F\\
\mathcal V/H & \xrightarrow{\hspace{0.5cm} \dot{k} \hspace{0.5cm}} & \mathcal V_0/H\times F
\end{array} \qquad \alpha_0 \circ k = \dot{k}  \circ \alpha, 
\]
where 
	$\alpha_0\times id_F(l_0, t)= (H.l_0,t) $  and $\dot{k}  (H.l)= (p_{\g, \g_0}(H.l), l.(X_0+ X_l))$.
	
	Therefore, 
	 the push-forward of $\dot\nu_{f, \h}^H$ under the map $\dot{k}$ is equivalent to the product measure
	$\dot\nu_{f_0, \h}^H\times d_Ft$.
	
	 Let $l_0 \in   \mathcal V_0$.
	For each $t \in F$, let  $l_{(t)}\in \g^*$ be defined   by 
		$$      p_{\g, \g_0}(l_{(t)})= l_0 \mbox{ and } l_{(t)}(X_0 + X_{{l_0}})= t.       $$
		In other words $ l_{(t)} = k^{-1}(l_0, t)$.
		By Theorem \ref{th1}, we have  
		$$      \operatorname{Ind}_{G_0}^{G}\pi_{l_0} =  \int^\oplus_F\pi_{k^{-1}(l_0, t)}  d_Ft .    $$
		So that 
		\begin{eqnarray*}
		\int^\oplus_{X^0_\chi/H}\operatorname{Ind}_{{G_0}}^{G} \pi_{l_0}  \diff\dot\nu_{f_0, \h}^H(H.l_0)   &=& 
		\int^\oplus_{X^0_\chi/H}\int^\oplus_F\pi_{k^{-1}(l_0, t)}  d_Ft \diff\dot\nu_{f_0, \h}^H(H.l_0) \\
		&=&  \int^\oplus_{X_\chi/H} \pi_l  \diff\dot\nu_{f, \h}^H(H.l).
		\end{eqnarray*}				
This proves the result.
\subsection{ Case II: Suppose that  $\Xi \cap \mathcal U = \varnothing$.}\label{ahm2}

The proof proceeds exactly as in Lipsman \cite[pp. 441–444]{lipsman-1989}.

\bigskip 

\begin{center}
{\bf Part II : The Plancherel Formula for  $  \operatorname{Ind}_{H}^{G}\chi$ }
\end{center}
\section{Distribution vectors : }\label{dist}

Let $(\pi, {\mathcal H_\pi})$ be a unitary representation of $G$. Denote by  $ \mathcal H_\pi^\infty$ the  $G$-invariant subspace of   all smooth vectors of $\pi$. The action of  $G$  on $ \mathcal H_\pi^\infty$  is denoted by $\pi_\infty$  (See section \ref{repst}).
 
In fact, let    $(K_n)_n$  be   a  decreasing sequence of  open compact  subgroups in $G$ such that $\displaystyle \bigcap_{n}K_n = \{1\}$.  Then we have $$ \displaystyle{\mathcal H}_{\pi}^{\infty} = \displaystyle\bigcup_{n} {\mathcal H}_{\pi}^{K_n},$$ where 
$   {\mathcal H}_{\pi}^{K_n}=\left\{ v \in {\mathcal H}_{\pi} \, \vert \, \pi(k) v =v , \,\forall k \in K_n \right\}$.
We equip $\mathcal H_{\pi}^{\infty}$ with its natural locally convex inductive-limit topology, so that $\mathcal H_{\pi}^{\infty}$ is an LF-space :
A subset $W$ of ${\mathcal H}_{\pi}^{\infty}$ is a neighborhood of zero in  ${\mathcal H}_{\pi}^{\infty}$ if and only if, for every $n$, $W\cap   {\mathcal H}_{\pi}^{K_n}$ is a neighborhood of zero in $ {\mathcal H}_{\pi}^{K_n}$ (See \cite[Section 3]{dieudonne-1949}). For this topology,  a linear map $u : {\mathcal H}_{\pi}^{\infty}\longrightarrow  \cit $ 
is continuous if and only if, for every $n$, the restriction  $u_{\vert  {\mathcal H}_{\pi}^{K_n}}$ of $u$ to $ {\mathcal H}_{\pi}^{K_n}$ is continuous. Moreover, the representation  $\pi_\infty$ is a continuous representation of $G$. Indeed, 
let $v \in  {\mathcal H}_{\pi}^{\infty}$.  There exists $n$ such that $v \in {\mathcal H}_{\pi}^{K_n}$. Since the topology on ${\mathcal H}_{\pi}^{K_n}$ is the same as the induced topology of ${\mathcal H}_{\pi}^{\infty}$ and that  ${\mathcal H}_{\pi}^{K_n}$ is a closed subset of ${\mathcal H}_{\pi}$, the mapping: $$ x \in G \longmapsto\pi(x)v\in  {\mathcal H}_{\pi}^{\infty},$$ is continuous at $1$, hence it is  continuous everywhere on $G$. 

 The topological anti-linear dual    of ${\mathcal H}^{\infty}_{\pi}$ is called 
 the space of distribution vectors of $\pi$ and is denoted by ${\mathcal H}_{\pi}^{-\infty}$.
 Here ${\mathcal H}_{\pi}^{-\infty}$ is equipped by the strong dual topology. By transposition, we define 
$ \pi_{-\infty}$, a   representation  of $G$ in ${\mathcal H}_{\pi}^{-\infty}$, by the following formula:
$$ \langle \pi_{-\infty}(g)a,v \rangle = \langle a, \pi_{\infty}(g^{-1})v \rangle,  \, \,   a \in
 {\mathcal H}_{\pi}^{-\infty}, \,  v \in {\mathcal H}_{\pi}^{\infty}   , \,   g \in G.$$
\begin{remm}
\begin{enumerate}
	\item 
 If $\pi$ is irreducible then $\pi$ is admissible so that, in this case, every linear form on ${\mathcal H}_{\pi}^{\infty}$ is continuous.
\item \label{anti} If we denote by  ${\mathcal H}_{\pi, L}^{-\infty}$  to be  the topological linear dual    of ${\mathcal H}^{\infty}_{\pi}$   equipped with  the strong dual topology.  Then the map
\begin{eqnarray}\label{antil}
  \Theta : {\mathcal H}_{\pi, L}^{-\infty}\longrightarrow {\mathcal H}_{\pi}^{-\infty}, && \, a \longmapsto\overline{a}
	\end{eqnarray}
is an anti-linear isomorphism. Here 
$  \overline{a}$ is defined by 
$$    \langle \overline{a}, v \rangle =  \overline{\langle {a}, v \rangle}, \, \, v \in  {\mathcal H}_{\pi}^{\infty}.       $$
\item In the real case, the space ${\mathcal H}_{\pi}^{\infty}$ is a Fr\'echet space.
\end{enumerate}
\end{remm}

\subsection{}\label{distr}  Fix a    Haar measure $\diff_G g$ on $G$. 
Let $\var \in \mathcal  C_c^\infty(G)$ and $T \in {\mathcal H}^{\infty}_{\pi}$.  

For every $v \in {\mathcal H}_{\pi}^{\infty}$, we have 
$$ \int_G \var(g) \langle \pi_{-\infty}(g)T, v \rangle   \diff_G g =   \langle T,  \int_G \overline{\var(g)}\pi_{\infty}(g^{-1}) v \rangle   \diff_G g  \rangle = \langle T, \pi_\infty(\check{\var}) v \rangle,    $$
where $\check{\var}(g)= \overline{\var(g^{-1})}$. 
 Let $n$ be an integer large enough such that $ \check{\var}$ is $K_n$-biinvariant. It follows that 
$$   \check{\var}= \sum_{k=1}^r c_k1_{x_k K_n}  \mbox{ and }  \check{\var}=e_{K_n}*\check{\var}* e_{K_n} ,       $$
where $ e_{K_n}= \frac{1}{\vol(K_n) } 1_{K_n}$. Then 
 $$ \pi_\infty(\check{\var}) v = \sum_{k=1}^r c_k \pi_\infty( e_{K_n}) \pi_\infty(x_k)\pi_\infty( e_{K_n})v \in  {\mathcal H}^{K_n}_{\pi}.$$  
On the other hand, the anti-linear form 
$$     {\mathcal H}_{\pi}^{K_n} \ni w \longmapsto  \langle T, w \rangle    $$
is   continuous. Therefore, there exists a unique vector $w_n \in {\mathcal H}_{\pi}^{K_n} $
such that 
\begin{equation*}
\langle T, w \rangle = \langle w_n, w       \rangle_{\mathcal H_\pi} , \, w \in {\mathcal H}_{\pi}^{K_n}.
\end{equation*}
Thus, for  every  $v \in {\mathcal H}_{\pi}^{\infty}$,
\begin{equation}\label{riezz}
  \langle T, \pi_\infty(\check{\var}) v \rangle =  \sum_{k=1}^r \overline{ c_k } \langle w_n,    \pi_\infty( e_{K_n}) \pi_\infty(x_k)\pi_\infty( e_{K_n})v    \rangle_{\mathcal H_\pi}     
	\end{equation}
Hence, the anti-linear form ${\mathcal H}_{\pi}^{\infty} \ni v \longmapsto  \langle T, \pi_\infty(\check{\var}) v \rangle    $
is  continuous. 

Throughout, the distribution vector  defined by  formula \ref{riezz} is denoted by $\pi_{-\infty}(\var)T$.
\begin{lem}\label{dis_vec} Let $\var \in  \mathcal C_c^\infty(G)$ and $T \in {\mathcal H}^{-\infty}_{\pi}$. 
Then the distribution vector $\pi_{-\infty}(\var)T$ belongs to ${\mathcal H}^{\infty}_{\pi}$, that is, there exists a unique vector $w \in {\mathcal H}_{\pi}^{\infty}$, such that 
$$      \langle \pi_{-\infty}(\var) T,  v \rangle = \langle w, v       \rangle_{\mathcal H_\pi} , \, \forall v \in {\mathcal H}_{\pi}^{\infty} .       $$
\end{lem}
\begin{dem} The result follows immediately from \eqref{riezz}.
\end{dem}
\begin{lem}
 $\mathcal H_\pi$  is dense in $ {\mathcal H}^{-\infty}_{\pi}$.
\end{lem}
\begin{Proof} Let $T \in  {\mathcal H}^{-\infty}_{\pi}$,  $B$ a bounded subset of ${\mathcal H}_{\pi}^{\infty}$,  and $\e>0$.

\noindent A standard property of strict LF-spaces is that every bounded subset is contained and bounded in one of the defining Fréchet spaces (See \cite[Proposition 4]{dieudonne-1949}). Then there exists an integer $n$ large enough such that 
$    B \subset {\mathcal H}_{\pi}^{K_n}$ and $B$ is bounded in  that space.  Let $w_n $ be the unique vector in  $ {\mathcal H}_{\pi}^{K_n} $
such that 
$$
\langle T,  v \rangle = \langle w_n, v       \rangle_{\mathcal H} , \, \forall v \in {\mathcal H}_{\pi}^{K_n}.
$$
Denote by $ T_{w_n} $ the distribution vector defined by $w_n$, that is, 
$$  \langle  T_{w_n},  v \rangle = \langle w_n, v       \rangle_{\mathcal H} , \, \forall v \in {\mathcal H}_{\pi}^{\infty}.      $$
We have $$ \sup_{v \in B}|\langle T-  T_{w_n}, v \rangle|  =0.           $$
 Thus the result.
\end{Proof}
{\bf Example :} Consider  the unitary representation $ \pi= \operatorname{Ind}_{H}^{G}\chi$ acting by  left translations   on the Hilbert space   $\mathcal H_\pi : =\mathcal H_{f,\h}$ . Let $a_\pi$ be the anti-linear form on $\mathcal H_{\pi}^\infty$ defined by :
\begin{equation}\label{dis_e}
a_\pi : \mathcal H_\pi^\infty \longrightarrow \cit, \, \, \var \longmapsto \overline{\var(e)}
\end{equation}
\begin{lem}\label{dist_e}  
\begin{enumerate}
	\item 
The anti-linear form $a_\pi$ belongs to $ {\mathcal H}^{-\infty}_{\pi}$ and satisfies 
$$  \pi_{-\infty}(h) a_\pi= \chi_f(h)a_\pi, \, \,  h \in H.       $$
\item For every $\phi \in \mathcal C_c^\infty(G)$, we have 
\begin{eqnarray*}
  &&\pi_{-\infty}(\phi) a_\pi (x)=  \int_H \phi(xh) \chi_f(h)\diff_Hh , \, x \in G    \\
		&& \langle \pi_{-\infty}(\phi) a_\pi, a_\pi \rangle = \int_H \phi(h) \chi_f(h)\diff_Hh \, \, \,  \mbox{ (matrix coefficient of $a_\pi$) } 
	\end{eqnarray*}
	\item The distribution vector $a_\pi$ is cyclic for $\pi$.
\end{enumerate}
\end{lem}
\begin{Proof}
Let $n \in \N\setminus\{0\} $. We should show that $a_{\pi\vert \mathcal H_\pi^{K_n}}$ is bounded. 

\noindent Because $K_n$ is compact and open in $G$, its image in $G/H$, namely $K_nH/H$, is compact and open.
As $\diff_{G/H}\dot x$ is  a regular Borel measure on $  G/H $, it follows that  
$$   0 < \vol(K_nH/H)     <\infty.         $$
Now, we have 
$$  |\langle a_\pi, \var \rangle | = |\overline{\var(e)}| \leq \frac{ ||\var||_{\mathcal H_\pi}}{vol(K_nH/H)}, \,  \var \in  \mathcal H_\pi^{K_n}. $$
The semi-invariance of $a_\pi$ and the second assertion of the lemma  are   immediate.

For the last assertion, it suffices to prove that if $\psi \in  \mathcal H_\pi^{\infty} $ satisfies 
$    \langle  a_\pi, \pi_\infty(x) \psi \rangle=0 , \, \forall x \in G$ then $\psi= 0$.   

Let  $\psi \in  \mathcal H_\pi^{\infty} $ be such a vector. Then, for every $  \phi \in \mathcal C_c^\infty(G)$, we have $\langle a_\pi, \pi_\infty(\phi) \psi \rangle=0 $. We deduce that 
$$     \langle \var, \psi         \rangle_{\mathcal H_\pi} =  0, \, \, \forall \var \in   \mathcal C_c^\infty(G, H, \chi),       $$
where $\mathcal C_c^\infty(G, H, \chi)$  is the space of complex-valued functions $\var : G\longrightarrow \cit$ satisfying the following conditions :
\begin{itemize}
	\item $\var(xh) = \chi_f(h^{-1}) \var(x)$, \, for $ h\in H, \,  x \in G $ ;
	\item  $\exists K$ an open compact subgroup of $G$ such that $\var(x_0x) = \var(x), \, \forall x \in G, x_0 \in K$ ; 
	\item $Supp (\var)$ is compact modulo $H$.
\end{itemize}
Since the subspace $\mathcal C_c^\infty(G, H, \chi) $ is dense in $ \mathcal H_\pi$, it follows that  $\psi= 0$.
\end{Proof}

\begin{lem} Suppose that the representation $(\pi, {\mathcal H_\pi})$ is admissible.

The  strong dual  of $ {\mathcal H}^{-\infty}_{\pi}$ is canonically isomorphic to ${\mathcal H}_{\pi}^{\infty}$, that is, the map 
$$  {\mathcal H}_{\pi}^{\infty} \ni  v   \longmapsto J(v) \in  ({\mathcal H}^{-\infty}_{\pi})',       $$
where
$\langle J(v), T \rangle =   \overline{\langle T,  v\rangle}$, is an isomorphism. 
\end{lem}
\begin{Proof}
The result  is essentially a reflexivity statement for the strict LF-space $\mathcal H_\pi^\infty$ with 
finite-dimensional stages.
\end{Proof}
In what follows, we use the notation $\langle v,T\rangle$ for $\langle J(v),T\rangle$.
\subsection{The  irreducible  case  : } 
\

 In this section,  let $ l\in \g^*$ and   $\b$  a polarization at $l$.  Recall that   $\chi_l$  is the character of the algebraic subgroup  $B = \exp(\b)$  defined by : 
\begin{equation*}
     \chi_l(\exp(X)) = \psi_F(\langle l, X\rangle), \, X \in \b. 
		\end{equation*}
		Denote by $\pi_l = \pi_{l,\b}$ the  induced representation of  $G$ : 
$     \pi_{l,\b} = \mbox{Ind}_B^G \chi_l$  acting in the Hilbert space  $\mathcal H_l = \H_{l,\b}$.   
Remark that this space is defined using  a $G$-invariant Borel measure $\diff_{G/B} \dot x$ on $G/B$. 
From \cite{dijk-1978}, the irreducible representation  $\pi_l $  is admissible. 

Assume that $G.l \cap (f + \h^\perp) \neq \varnothing$. Fix  an element $g \in G$ such that $ g.l \in   f + \h^\perp $.

  Let  $\diff_{H/H \cap gBg^{-1}} \dot{h}$ be a positive    $H$-invariant measure on  
 $H/H\cap gBg^{-1}$.  From\cite[Proposition 4.6.1]{maktouf-2026}, the anti-linear form
\begin{equation}\label{khema}
  I_{g.l, \pi_l}:  {\mathcal H}_{l, \b}^{\infty} \ni \phi   \longmapsto \overline{\int_{H/H \cap gBg^{-1}} \phi(hg)\chi_f(h)\diff_{H/H \cap gBg^{-1}}\dot h}
\end{equation}
is  non-zero.  Moreover such distribution vector satisfies : 
$$        \pi_{l, -\infty}(h) I_{g.l, \pi_l} = \chi_f(h)I_{g.l, \pi_l}, \, \forall h \in H.     $$ 
For $\phi\in \mathcal C_c^\infty(G)$,  Lemma \ref{dis_vec} implies that $  \pi_{l, -\infty}(\phi) I_{g.l, \pi_l}   $ belongs to $  {\mathcal H}_{l, \b}^{\infty} $. We will compute this function explicitly.

To this end, let $\phi_H^f$ be the complex-valued function on $G$ defined by : 
$$     \phi_H^f(x)= \int_H \phi(xh) \chi_f(h) \diff_H  h  , \, x \in G    $$
where $\diff_H  h$  denoted a Haar measure on $H$.
Then $ \phi_H^f$ is   locally constant  on $G$, has   compact support modulo $H$, and 
 satisfies :
$$       \phi_H^f(xh ) = \chi_f(h^{-1})\phi_H^f(x) , \, \, x \in G, h \in H. $$
$\bullet $ Choice of measures : 
\begin{enumerate}
	\item To define $\phi_H^f$, we fixe a Haar measure $\diff_H  h$  on $H$. 
	\item  To define $I_{g.l, \pi_l}$, we  fixe    a positive  $H$-invariant measure $\diff_{H/H \cap gBg^{-1}} \dot{h}$  on  
 $H/H\cap gBg^{-1}$.

We deduce that there exists a Haar measure $\diff_{H \cap gBg^{-1}}h$ such that 
$$ \diff_H  h =   \diff_{H/H \cap gBg^{-1}} \dot{h}  \diff_{H \cap gBg^{-1}}h.     $$
\item To define $ \pi_{l}(\phi) $,  we fix  a Haar measure $\diff_G  x$  on $G$. 
Then there exists a Haar measure $\diff_B b$ on $B$ such that 
$$   \diff_G  x = \diff_{G/B}\dot x \diff_B b.         $$
\item Let $    \diff_{ g^{-1}Hg \cap B}h$ be the Haar measure on $ g^{-1}Hg \cap B$ image of $\diff_{H \cap gBg^{-1}}h         $ by the group isomorphism : 
$$   H \cap gBg^{-1}\longrightarrow   g^{-1}Hg \cap B, \, h \longmapsto   g^{-1}hg  .    $$
\item We  equip  $   B/B\cap g^{-1}Hg    $ with the $B$-invariant measure $\diff_{B/B\cap g^{-1}Hg }\dot b$ such that 
$$    \diff_B b =   \diff_{B/B\cap g^{-1}Hg }\dot b    \diff_{ g^{-1}Hg \cap B}h.     $$
\end{enumerate}

\begin{pro}\label{liss} Let the notations be as above. For every $ \phi \in \mathcal C_c^\infty(G)$, we have 
$$   \pi_{l, -\infty}(\phi) I_{g.l, \pi_l} (x)= \int_{B/B\cap g^{-1}Hg} \phi_H^f(xbg^{-1}) \chi_l(b) \diff_{B/B\cap g^{-1}Hg} \dot b, \, x \in G $$
Moreover 
\begin{eqnarray}\label{beau}
&&\langle   \pi_{l, -\infty}(\phi) I_{g.l, \pi_l},   I_{g.l, \pi_l}  \rangle \\
&&\qquad   = 
  \int_{H/H \cap gBg^{-1}}\left(\chi_f(h)  \int_{B/B\cap g^{-1}Hg}   \phi_H^f(hgbg^{-1})\chi_l(b) \diff_{B/B\cap g^{-1}Hg} \dot b  \right)      \diff_{H/H \cap gBg^{-1}}\dot h    \nonumber
\end{eqnarray}
\end{pro}
\begin{Proof}
The result follows by a straightforward computation.
\end{Proof}

\section{Plancherel Formula}{}\label{Pl}
\

Throughout the remained of this paper,   we assume  that $H$ is a normal algebraic subgroup of $G$. 
Let $f$ be a linear form on $\g$ such that $\h$ is subordinate to $f$.  We consider the character 
$ \chi = \chi_f$ associated to $f  $,  and  the induced representation $ \pi=  \operatorname{Ind}_{H}^{G}\chi$  acting in the Hilbert space  $\H_{l,\h}$.

To this end, we first   investigate the $H$-orbits in  $f + \h^\perp$. Let $l \in f + \h^\perp$. Since $\h$ is an ideal of $\g$, we have $X.f =   X.l     $, for every $ X \in \h$. We deduce that $\h(l)= \h(f)$. 
On the other hand, since $\h$ is subordinate to $l$, we have 
$$ \exp(X).l = l+ X.l = l + X.f, \, \forall X \in \h.             $$
It follows that 
$$  H.l = l+ \h.f.        $$
In particular, all $H$-orbits in $f + \h^\perp$ have the same dimension, say $d$. 

By Lemma \ref{dim}, we have 
$   \dim G.l - \dim H.l \geq   \dim H.l$ which is equivalent to 
$$ \dim (\g(l)) \leq \dim \g -2d.       $$
Moreover,   equality  holds if and only if $ l \in    \Xi  $. We deduce that 
$$ \Xi  =\{l \in f + \h^\perp \, \vert \, \dim (\g(l)) = \dim \g -2d \}$$  and is  closed  in $f + \h^\perp $.

Let     $\mathcal U$ be the  nonempty Zariski-open set  of $f + \h^\perp$  introduced in Proposition  \ref{Zariski}. 
 By Proposition \ref{Zariski},  the set  $ \mathcal U \cap \Xi$ is Zariski open in $f + \h^\perp$. Since   $f + \h^\perp$ is irreducible, we deduce  that  $  \mathcal U \cap \Xi  \neq \varnothing$ if and only if  $ \Xi=f + \h^\perp$. 
If that is the case,  all  the multiplicities appearing in formula \ref{prin} are  finite.
Moreover, we  proved in \cite[Section 3.2]{maktouf-2026} that,  for every $l \in f + \h^\perp = \Xi$,  $  (f + \h^\perp)\cap G.l = H.l $.  Consequently all the  multiplicities appearing in formula \ref{prin} are equal to one.
\subsection{}\label{inv} In what follows, we assume that $ \Xi=  f + \h^\perp$.

  Recall that a  distribution  vector  $a \in {\mathcal H}_{\pi}^{-\infty}$ is said to be  $\chi$-semi-invariant if it    satisfies :
$$  \pi_{-\infty}(h)a = \chi(h)a, \,  \forall h \in H.$$
The space of such distribution vectors will be denoted by   $( {\mathcal H}_{\pi}^{-\infty})^{H, \chi}$. 

Let $l \in  f + \h^\perp$. Set $\b_l = \h + \g(l)$, which is a polarization at $l$, and let $B_l = \exp(\b_l)$ the corresponding algebraic subgroup.  Denote by $\pi_l = \operatorname{Ind}_{B_l}^{G}\chi_l$ the induced representation of $G$ by the character $\chi_l$ and realized in the Hilbert space   $\mathcal H_{l,\b_l}$. 

Following \cite[Proposition 5.2.1]{maktouf-2026},  $( \mathcal H_{l,\b_l}^{-\infty})^{H, \chi}$ is one dimensional vector space spanned by
the   distribution vector $I_{l, \pi_l}$ given by formula \ref{khema}. We have
  $$I_{l, \pi_l}:  \phi \in  {\mathcal H}_{l,\b_l}^{+\infty}    \longmapsto \overline{\phi(e) }. $$
	We equip $\h(f)$ with  a Haar measure   $\diff_{\h(f)}X$ and  $\h/\h(f)$ with a Haar measure   $ \diff_{\h/\h(f)}\dot X   $ such that $\diff_{\h}X = \diff_{\h/\h(f)}\dot X \diff_{\h(f)}X$. We equip  $\h.f$ with the Haar measure 
	$\diff_{\h.f}l' $ obtained as the image of   $ \diff_{\h/\h(f)}\dot X   $ under  the isomorphism 
	$$ \h/\h(f) \ni \dot X \longmapsto X.f \in   \h.f.    $$
	Finally, we equip $H.l$ with the $H$-invariant measure $ \diff_{H.l}l'$,  which is the image of 
	$ \diff_{\h.f}l'  $ under the diffeomorphism 
	$$     \h.f \ni X.f  \longmapsto  l + X.f.   $$
	\begin{pro} We have 
	$$   \langle    \pi_{l, -\infty}(\phi) I_{l, \pi_l},   I_{l, \pi_l} \rangle  =   \int_{H.l}\widehat{\phi \circ\exp}_\g( l')\diff_{H.l}l' , \, \,  \phi \in \mathcal C_c^\infty(G)      $$
	
	\end{pro}
 \begin{Proof}
Let  $ \phi \in \mathcal C_c^\infty(G) $.  By Proposition \ref{beau}, one has,  
\begin{eqnarray*}
\langle    \pi_{l, -\infty}(\phi) I_{l, \pi_l},   I_{l, \pi_l} \rangle  & =& \int_{B/H} \phi_H^f(b)\chi_l(b)   \diff_{B/H} \dot b   = \int_B\phi(b)\chi_l(b)\diff_Bb  \\
&=& \int_{\b}\phi \circ\exp(X)\psi_F(\langle l , X \rangle )\diff_\b X  \\
&=&\widehat{(\phi \circ\exp)}_\b(l) = \int_{\b^\perp}\widehat{(\phi \circ\exp)}_\g(l + l')\diff_{\b^\perp}l'\\
&=& \int_{B.l}\widehat{(\phi \circ\exp)}_\g( l')\diff_{B.l}l'=\int_{H.l}\widehat{(\phi \circ\exp)}_\g( l')\diff_{H.l}l'
\end{eqnarray*}
Thus the result.
\end{Proof}
\subsection{} \label{planch}
Let $W$ be a complementary subspace to $\h.f$ in $\h^\perp$. Then the map 
$$\Theta :  W\times H.f      \longrightarrow   f + \h^\perp, (w,l)\longmapsto w+l   $$
is a diffeomorphism.

 Let $\diff_Ww$ be the Haar measure on $W$ such that
$    \diff_Ww  \diff_{H.f}l' = \mu_{f, \h^\perp} $.   
\begin{tho} For every  $ \phi \in \mathcal C_c^\infty(G) $, we have 
$$   \langle \pi_{-\infty}(\phi) a_\pi, a_\pi \rangle  = \int_W    \langle    \pi_{w +f, -\infty}(\phi) I_{w+f, \pi_{w+f}},   I_{w+f, \pi_{w+f}} \rangle    \diff_Ww  .     $$
\end{tho}
\begin{Proof} Let  $ \phi \in \mathcal C_c^\infty(G) $. 
From the above computation, we have 
\begin{eqnarray*}
 \int_W \langle    \pi_{w +f, -\infty}(\phi) I_{w +f, \pi_{w +f}},  I_{w +f, \pi_{w +f}} \rangle \diff_Ww & =&   \int_W  \int_{H.f}\widehat{(\phi \circ\exp)}_\g( w+l)\diff_{H.f}l  \diff_Ww\\
&=& \int_{ f + \h^\perp} \widehat{(\phi \circ\exp)}_\g( l)  \diff\mu_{f, \h^\perp}(l) \\
&=& \int_{  \h^\perp} \widehat{(\phi \circ\exp)}_\g( f+l)  \diff_{ \h^\perp}(l)\\
&=&\int_{  \h} \phi \circ\exp(X) \chi_F(\langle f, X \rangle) \diff_{ \h}(X) \\
&=&  \langle \pi_{-\infty}(\phi) a_\pi, a_\pi \rangle.
\end{eqnarray*}
The desired conclusion follows.
\end{Proof}

\section{Appendix A :}{\bf Recall of  results of Pukanszky \cite[Propositions 1.1 and 4.1]{pukanszky-1971}}\label{apA}

 Let $G$  be the set of rational points   of a unipotent algebraic group defined over     $F$  and $\g$ be its  Lie algebra. The group $G$ acts on $\g^*$ via the coadjoint representation. 
\begin{tho}\label{AP}
f
Let $\g^* = V_0 \supset  V_1 \supset \ldots \supset V_m= 0  $ be a Jordan-H$\ddot{\mbox{o}}$lder sequence for $G$ and assume that $v_j \in V_{j-1} \setminus V_j       $, $1 \leq j \leq m$. Then there is a sequence of subsets $\g^* = F_0 \supset F_1   \supset \ldots \supset F_M \neq 0$ with the following properties : 
\begin{itemize}
	\item For any $j \,  (1 \leq j \leq M)$, $F_j$ is $G$-invariant and the dimension of any $G$-orbit in $ F_{j-1} \setminus F_j $ is constant ; $F_M$ is the collection of all fixed points of $G$.
	\item There is a collection of homogeneous polynomials $\{Q_j(x) ; 1 \leq j \leq M\}      $ on $\g^*$, such that $ F_j= \{x ; x \in\g^* , Q_k(x)= 0 \mbox{ for } k \leq j \} $.
	\item Let $j$ be fixed as above, and let us denote by $d$ the dimension of some  $G$-orbit in $ F_{j-1} \setminus F_j $. There is a subset $ 0<j_1 < \ldots<j_d<m           $ of $\{1,2, \ldots, m\} $ and a system of $m$ functions $\{P_k(z,x) ; 1 \leq k \leq m\}$ on $ F^d \times  ( F_{j-1} \setminus F_j)  $ with the following properties : 
	\begin{enumerate}
		\item For sufficiently large integer $N$, $(Q_j(x))^NP_k(z,x) \, (1 \leq k \leq m)$ is the restriction
		of polynomial function on $ F^d \times \g^*  $  to  $ F^d \times  ( F_{j-1} \setminus F_j)  $. 
		\item \label{a11} We have $ P_{j_k}(z,x)= z_k $ , $k=1, \ldots, d$.
		\item \label{b11}  If $i$ is some integer between $1$ and $m$ and $t = \sup_{j_k \leq i}k$, $P_i(z,x)$ depends on $\{z_h ; h \leq t\}     $ only.
		\item  \label{c11} For any $x \in  F_{j-1} \setminus F_j$ we have 
		$$G.x = \{v  ; v = \sum_{k=1}^m P_k(z,x).v_k, \, z \in F^d\}.           $$
		For each fixed  $x$, the functions  $\{P_k(z,x) ; 1 \leq k \leq m\}$ are uniquely determined by
	\end{enumerate} the condition, that they are polynomials in $z \in F^d$ satisfying conditions \eqref{a11}, \eqref{b11} and \eqref{c11}.
	\end{itemize}
 \end{tho}

\begin{co}\label{AC} Keep the notations of Theorem \ref{AP}.

The map $\lambda :  F_{j-1} \setminus F_j\longrightarrow \g^*$ defined by 
$$ \lambda(x)= \sum_{k=1}^m  P_k(0,x) v_k              $$
is regular, $G$-invariant, and  verifying $  \lambda(x) \in G.x   $.  
\end{co}

\bibliographystyle{amsplain} \bibliography{biblio}

\end{document}